\documentclass[12pt,a4paper,twoside,final,notitlepage, leqno]{article}
\usepackage[english]{babel}
\usepackage[T1]{fontenc}  
\usepackage{epsfig, graphicx, amssymb}
\usepackage{amsmath,amsthm,epsfig,amsfonts}  
\usepackage{float}
\usepackage{color}
\def\br {\break}

\newcommand{\monitem}{ \smallskip \noindent $\bullet$ \quad  } 
\newcommand{\moneq}{\vspace*{-7pt} \begin{equation} \displaystyle } 
\newcommand{\moneqstar}{\vspace*{-6pt} \begin{equation*} \displaystyle } 
\newcommand{\monendstar}{\vspace*{-6pt} \end{equation*}   }
\newcommand{\monend}{\vspace*{-7pt} \end{equation}   }

\newcommand{\dd}{{\rm d}}

\newcommand{\RR}[0]{\mathbb{R}}

\def\section*#1{}

\usepackage{fancyhdr}
\renewcommand{\headrulewidth}{0pt}

\begin{document} 

\fancypagestyle{plain}{ \fancyfoot{} \renewcommand{\footrulewidth}{0pt}}
\fancypagestyle{plain}{ \fancyhead{} \renewcommand{\headrulewidth}{0pt}}

~

  \vskip 2.1 cm

\centerline {\bf \LARGE  Unexpected convergence  }

 \bigskip 

\centerline {\bf \LARGE  of lattice Boltzmann schemes }

 \bigskip  \bigskip \bigskip

\centerline { \large    Bruce M. Boghosian$^{a}$, Fran\c{c}ois Dubois$^{bc}$, Benjamin Graille$^{b}$,}

\smallskip  \centerline { \large   Pierre Lallemand$^{d}$ and Mohamed-Mahdi Tekitek$^{e}$}

\smallskip  \bigskip 

\centerline { \it  \small   
$^a$  Dpt. of Mathematics, Tufts University, Bromfield-Pearson Hall, Medford, MA 02155, U.S.A.} 

\centerline { \it  \small   
$^b$   Dpt. of Mathematics, University Paris-Sud,  B\^at. 425, F-91405  Orsay, France.} 

\centerline { \it  \small   
$^c$    Conservatoire National des Arts et M\'etiers, LMSSC laboratory,  F-75003 Paris, France.} 

\centerline { \it  \small  $^d$   Beijing Computational Science Research Center, 
Haidian District, Beijing 100094,  China.}

\centerline { \it  \small  $^e$   Dpt. Mathematics, Faculty of Sciences of Tunis, 
University Tunis El Manar, Tunis, Tunisia. } 


\bigskip  

\centerline { 30 April 2018 
{\footnote {\rm  \small $\,$ Contribution  presented to the 
    26th DSFD Conference, Erlangen (Germany), 10 - 14 July 2017.
  {\it Computers \& Fluids}, volume 172,   pages 301-311,  2018.}}}


 
 {\bf AMS (MSC2010) classification}: 76M28.

\bigskip  
\noindent {\bf \large Abstract} 

\noindent 
In this work, we study 
numerically 
the convergence of the scalar D2Q9 lattice Boltzmann scheme 
with multiple relaxation times when the time step  is proportional to the space step and tends to zero. 
We do this by a combination of theory and numerical experiment.  
The classical formal analysis when all the relaxation parameters are fixed 
and the time step tends to zero shows 
that the numerical solution converges to solutions of the heat equation, 
with a constraint connecting the diffusivity, the space step and the coefficient of relaxation 
of the momentum.
If the diffusivity is fixed and the space step tends to zero, 
the relaxation parameter for the momentum is very small, 
causing a discrepency between the previous analysis and the  numerical results.  
We propose a new analysis of the method for this specific situation of evanescent relaxation, 
based on the dispersion equation of the lattice Boltzmann scheme. 
 A new asymptotic partial differential equation, the damped acoustic system, 
is emergent as a result of this formal analysis.  
Complementary numerical experiments establish the 
convergence of the scalar D2Q9 lattice Boltzmann scheme 
with multiple relaxation times and acoustic scaling in this specific case 
of evanescent relaxation towards the numerical solution of the damped acoustic system.

\bigskip \bigskip   \newpage \noindent {\bf \large    1) \quad  Introduction  }    

\fancyhead[EC]{\sc{ B.M. Boghosian, F. Dubois,  B. Graille, P. Lallemand, and  M.M. Tekitek }} 
\fancyhead[OC]{\sc{Unexpected convergence  of lattice Boltzmann schemes }} 
\fancyfoot[C]{\oldstylenums{\thepage}}

\noindent 
Lattice Boltzmann models are simplifications of the continuum Boltzmann 
equation obtained by discretizing in both physical space and velocity space.  
The discrete velocities $v_i$ retained typically correspond to lattice vectors 
of the discrete spatial lattice.  That is, each lattice vertex $x$ is linked 
to a finite number of neighboring vertices by lattice vectors $v_i \, \Delta t$.
A particle distribution $ \, f \, $ is therefore parametrized by its components 
in each of the discrete velocities, the vertex $x$ of the spatial lattice, 
and the discrete time $t$.  A time step of a classical lattice Boltzmann 
scheme~\cite{LL00} then contains two steps:

\noindent
{\it (i)} a relaxation step where the distribution $ f $ at each vertex $x$ 
is locally modified into a new  distribution $ f^* $, and

\noindent
 {\it (ii)} an advection step based  on the method of characteristics 
as an  exact time-integration operator.  We employ the  multiple-relaxation-time 
approach introduced by d'Humi\`eres~\cite{DDH92}, wherein 
the local mapping $ \, f \longmapsto  f^*  \, $ is  described by a  
diagonal  operator in a space of moments. 

\noindent 
In~\cite{DL09}, we have 
studied the asymptotic expansion of various 
lattice Boltzmann schemes with multiple-relaxation times 
for different applications.  We used the so-called acoustic scaling, in which
the ratio  $ \, \lambda \equiv \Delta x / \Delta t \,$ is kept fixed. 
%
We supposed also that the relaxation operator remains fixed. 
%
In this manner, we demonstrated
the possibility of approximating diffusion processes described by the heat equation.

\noindent 
The importance of using small values of relaxation parameters
was 
recognized for  linear viscoelastic fluids by Lallemand {\it et al.} 
\cite {LHLR03}.
Independently, unexpected results 
in simulations  for advection-diffusion processes have been  described by 
Dellacherie in~\cite{De14}. 
We have studied experimentally in~\cite {BDGLT18} 
the curious convergence of the D1Q3 
multiple-relaxation time lattice Boltzmann scheme 
with one conserved variable when using the  acoustic scaling 
%
in one spatial dimension. 
The asymptotic equation of the lattice Boltzmann scheme
is no longer an advection-diffusion model but a damped acoustic model. 
In this contribution, we show and analyze an analogous phenomenon for  two 
spatial   
 dimensions with the scalar D2Q9 lattice Boltzmann scheme.
The difficulty concerns the  highlighting of the convergence with the numerical experiments. 
%

\noindent 
In Section~2, we recall some fundamentals relative to the D2Q9 lattice
Boltzmann scheme for scalar conservation laws. In Section~3, we study convergence 
of this scheme for diffusive and acoustic scaling. A formal analysis 
is proposed in Section~4, with the dispersion equation method, 
initially proposed in \cite{Pi68}.
We establish that with acoustic scaling, the convergence of the 
scalar D2Q9 scheme is not the heat equation but an unexpected model! 
Finally, we study  the experimental  convergence of the scalar D2Q9 
scheme in several 
situations in Section 5.

 \bigskip \bigskip   \noindent {\bf \large    2) \quad  
Scalar  D2Q9 lattice Boltzmann scheme for thermal problems }   

\begin{figure}    [H]  \centering 
\includegraphics[width=.45 \textwidth, angle=0]{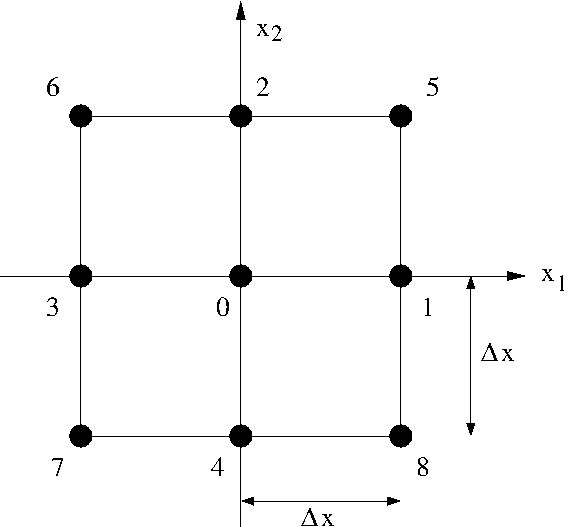} 
\caption{Particle distribution $ \, f_j \, $ for $ \, 0 \leq j \leq 8 \, $
of the D2Q9 lattice Boltzmann scheme} 
\label{d2q9stencil} \end{figure}

\noindent  
The D2Q9 lattice Boltzmann scheme  uses a set of discrete velocities described 
in Figure~\ref{d2q9stencil}. 
A density distribution $ \, f_j \,$ is associated to each velocity
 $ \, v_j  \equiv \lambda \, e_j $,  
where 
 $ \, \lambda \equiv {{\Delta x}\over{\Delta t}} \, $ is the  fixed numerical lattice velocity.
The first three moments for the 
density and momentum are defined according to 
\moneq \label{moments-conserves} 
\rho \, = \,  \sum^8_{j=0} f_{j} \, = \, m_0 \,, \,\,\,  
J_x \equiv \rho \, u_x \,  = \,   \sum^8_{j=0} \lambda \, e^1_{j} \, f_{j} \, = \, m_1 \,, \,\,\,    
J_y \equiv \rho \, u_y \, = \,   \sum^8_{j=0} \lambda \, e^2_{j} \, f_{j} \, = \, m_2 \, ,  
\monend    
where the  $\, e^\alpha_{j} \, $ are the  
$\alpha$th        
cartesian components  
of the vectors $e_{j}$ introduced previously. 
We complete this set of moments and 
construct a vector  $ \, m \, $  of moments  $ \, m_k \,$ according to 
\moneq \label{moments} 
m = M \, f \,, 
\monend
with an invertible fixed  matrix $ \, M \,$ usually \cite{LL00}  given by 
\moneqstar 
 M  \, = \, \left (\begin {array}{ccccccccc}
\displaystyle   1 &   1 & 1 & 1 & 1 & 1 & 1 & 1 & 1\cr
\displaystyle   0 &   \lambda & 0 & -\lambda & 0 & \lambda & -\lambda & -\lambda  & \lambda\cr
\displaystyle   0 &   0 & \lambda & 0 & -\lambda & \lambda & \lambda & -\lambda & -\lambda\cr
\displaystyle -  4 \, \lambda^2 &   -\lambda^2  &  -\lambda^2 &  -\lambda^2 & -\lambda^2 &
 2 \,\lambda^2  &  2 \,\lambda^2 &  2 \,\lambda^2 &  2 \,\lambda^2 \cr  
\displaystyle   0 &    \lambda^2 & -  \lambda^2 &  \lambda^2 & -  \lambda^2 & 0 & 0 & 0 & 0\cr
\displaystyle   0 &   0 & 0 & 0 & 0 & 
\lambda^2 & -\lambda^2 & \lambda^2 & -\lambda^2    \cr  
\displaystyle   0 &   -2 \, \lambda^3 & 0 & 2  \, \lambda^3 & 0 & \lambda^3 &
 - \lambda^3 & - \lambda^3 &  \lambda^3 \cr
\displaystyle   0 &   0 & -2 \, \lambda^3 & 0 & 2 \, \lambda^3 &
  \lambda^3 &  \lambda^3 & - \lambda^3 & -  \lambda^3 \cr
\displaystyle   4 \, \lambda^4 &  -2 \, \lambda^4 &  -2  \, \lambda^4& -2  \, \lambda^4 &
 -2  \, \lambda^4 &  \lambda^4 & \lambda^4  &  \lambda^4 &   \lambda^4   \end{array}  \right) \, . 
\monendstar 
For scalar lattice Boltzmann  applications, the density $ \, \rho \, $
is the   ``conserved variable''. 

\monitem 
The particle distribution at equilibrium  $ \, f^{\rm eq} \, $ is a function only
of this conserved variable. For this thermal D2Q9 lattice Boltzmann scheme, 
the vector of equilibrium moments $ \, m^{\rm eq} \, $ 
is given by
\moneq \label{moments-equilibre} 
m^{\rm eq}  \,= \, \big( \rho \,,\, 0 \,,\, 0 \,, \, \alpha \, \lambda^2 \, \rho  
 \,,\, 0 \,,\, 0 \,, \, 0 \,,\, 0   \,,\,  \lambda^4 \, \beta \, \rho \big)^{\displaystyle \rm t} \, . 
\monend    
In most  applications, the coefficients $ \, \alpha \, $ and $ \, \beta \, $ are
usually taken  to be
\moneq \label{alpha-beta} 
\alpha = -2 \,,\,\,\, \beta = 1  \, . 
\monend    
The lattice Boltzmann scheme is comprised of  two  fundamental steps~: relaxation and advection.
During the relaxation step, the conserved variable  $ \, \rho  \, $ 
is   not   modified,  and the  non-conserved  moments $\, m_1 \, $ to  $\, m_8 \, $ 
relax towards  an equilibrium value: 
 $ \,\, \displaystyle   m_k^{\rm eq} \, = \, \psi_k (\rho) \,\, $ for $ \, k  \geq 1 $, 
where the  $\psi_k$ are the linear functions 
of the conserved moment given by (\ref{moments-equilibre}).
The specification of this  step also needs  relaxation rates 
$ \, s_k \, $ : for $ \, k \geq 1 \, $ such that 
\moneqstar 
 m_k^* \, = \,  m_k \,+\, s_k \, \big(  m_k^{\rm eq} \,-\,  m_k \big)  \,, 
\monendstar 
where the superscript $*$ denotes the moment $m_k$ after the  relaxation step.    
The table of relaxation parameters $ \, s_k \, $ chosen in our simulations is as follows 
\moneq \label{coef-relax} 
[s]  \,= \, \big( s_J \,,\, s_J \,,\, s_e \,,\, s_x \,,\, s_x \,,\, 
s_q \,,\, s_q \,,\, s_\varepsilon \big)  \, . 
\monend     
We introduce also the $8\times 8$  diagonal matrix $ \, S \, $ whose diagonal elements are the components of
 the vector $ \, [s] $. 
In our computations, we take  the following numerical values  
\moneq \label{coef-relax-juillet2017} 
s_e = 1.7 \,,\,\, s_x = 1.1 \,,\,\, s_q = 1.1 \,,\,\,   s_\varepsilon = 1.7 \, . 
\monend   
%
Only the relaxation coefficient $ \, s_J  \, $ for the first order momentum
is allowed to vary  in our numerical  experiments. 

\noindent
Then using the matrix $M^{-1}$ the relaxation step becomes in $f$ space~:
\moneq \label{from-m-to-f}
 f^*_i(x, \, t) \, = \, \sum_\ell \, M^{-1} _{i \, \ell} \,  m_\ell^{*}.
\monend 
During the advection step  $f_i(x_j)$ is transported from the node $x_j$    
by the discrete velocity $v_i$ to the node $x_j+v_i \Delta t.$
Thus the evolution of populations $f_i$ for $0\leq i \leq 8$  at internal node  
$x$ is described by:
\moneq \label{LB-schema} 
f_i(x,t+\Delta t) \, = \, f^*_i(x - v_i\Delta t,t) \, , \,\, 0 \leq i \leq 8 \, .
\monend

\monitem
In~\cite{DL09},
we have analyzed several lattice Boltzmann models 
with the Taylor-expansion  method, including the present one defined by
Eqs.~(\ref{moments}, \ref{moments-equilibre}, \ref{coef-relax},  \ref{LB-schema}).
The hypothesis used was that the reference  velocity $ \, \lambda \,$ 
and the relaxation coefficients $ \, s_J $, 
$ \, s_e $, $ \,  s_x $, $\, s_q \, $ and $ \, s_\varepsilon \, $ remain constant 
as the spatial step $ \, \Delta x \, $ tends to zero. 
Then the  conserved variable~$ \, \rho \, $  satisfies (at least formally!)
the heat equation:
\moneq \label{edp-diff} 
{{\partial \rho}\over{\partial t}}   \, - \, \kappa \,\,  
\Delta \rho   \,=\, {\rm O}(\Delta x^2)  \, ,
\monend
where the thermal diffusivity  $ \kappa  \, $ is given by the relation  
\moneq \label{mu} 
\kappa  \,\equiv\, {{4+\alpha}\over{6}} \, \sigma  \, \lambda \,  \Delta x \,,\quad  
\sigma \, \equiv \, \Big({{1}\over{s_J}} - {{1}\over{2}} \Big) \, . 
\monend
The coefficient $ \, \sigma \,$ is known as the ``H\'enon parameter''
in reference to the pioneering work of 
H\'enon~\cite{He87}. 
Observe that when the relaxation coefficient $ \, s_J \, $ and the mesh velocity
$ \, \lambda \, $ are fixed, the  thermal diffusivity tends to zero as the space step 
$ \, \Delta x \, $ tends to zero. 
%
This lattice Boltzmann scheme is 
stable  in the fluid case (see \cite{LL00}) under the condition: 
\moneqstar 
-4 \, < \, \alpha \, < \, 2 \, .  
\monendstar
For the scalar case, the condition $ \, \alpha + 4 > 0 \, $ is clear to assume
that the thermal diffusivity~$ \, \kappa \, $ is positive
(see (\ref{mu})) and the condition $ \, \alpha \, < \, 2 \, $ corresponds to our
experimental know how. 
%
Observe that with these choices, the value of the relaxation parameter  
$ \, s_J \, $ has to be fit with the physical diffusivity $ \, \kappa \, $
and the mesh size $ \, \Delta x \, $ through the relation (\ref{mu}) 
if the space step and time step are varying proportionately. 
In particular, we have the expansion   
\moneq \label{s-kappa} 
s_J \,=\, {{4+\alpha}\over{6 \, \kappa }} \,  \lambda \,  \Delta x \,+\, {\rm O}(\Delta x^2)  \,
\monend
as $ \, \Delta x \, $ tends to zero.

\monitem
Diffusive scaling can also be used and we refer, e.g.,  to the work of  
Junk {\it et al.} 
\cite{JKL05}. 
In this case, the ratio 
\moneqstar 
{\left({\Delta x}\right)^2\over{\Delta t}}
= \lambda \, \Delta x  
\monendstar
remains fixed. 
This diffusive scaling is intensively used with the explicit finite
difference method for solving the heat equation. 
It is well known \cite{RM57} that the time step must be
proportional to the square of the 
spacial 
step  in order for the method to be stable. 
An asymptotic analysis can be done for this simple lattice Boltzmann 
thermic model, as, e.g.,  in our contribution \cite{DL13},  
and we obtain again the heat equation (\ref{edp-diff}) as the
scaling limit of the model. 
With this diffusive scaling, the parameters $ \, \sigma \, $ and $ \, s_J \, $
remain constant 
%
if the  thermal diffusivity is given and the mesh size 
$ \, \Delta x \, $ tends to zero. 
Remark also that the convergence of the lattice Boltzmann scheme was rigorously proved for
the diffusive scaling for Navier-Stokes flows in periodic and bounded domains in \cite{JY09}
and for one  dimensional convection-diffusion-reaction equations in \cite{JY15}. 


\bigskip \bigskip   \noindent {\bf \large    3) \quad   First numerical experiments }   

\noindent 
%
We study the diffusion of a Gaussian profile in a square 
domain. 
In order to control
the computer cost during the numerical experiment and to be certain that the numerical experiment
is not polluted by the boundary scheme, we impose periodic boundary conditions.
We use  two variants of the scalar D2Q9 lattice Boltzmann scheme:
diffusive and  acoustic scaling. 

\monitem 
 Scalar D2Q9 numerical experiments with diffusive scaling  

\noindent
We solve numerically the heat equation 
\moneq \label{chaleur} 
{{\partial \rho}\over{\partial t}}   \, - \, \kappa \,\,  
\Delta \rho   \,=\, 0 \, ,
\monend
in the square $ \, \Omega = \, [ -1 , \, 1 ]^2 $, 
with periodic boundary conditions. 
The initial condition is a Gaussian:   
\moneq \label{gaussienne}  
\rho_0 (x,\, y) = \exp \Big(- {{x^2 + y^2}\over{0.09}} \Big) \,, \,\, -1 \leq x , \, y \leq 1 \, . 
\monend
The coefficients $ \, \alpha \, $ and $ \, \beta \, $ of the equilibrium are 
fixed according to (\ref{alpha-beta}) and  we keep fixed the 
relaxation coefficient for momentum :  
\moneq \label{sJ-fixe} 
s_J = {3\over2}  \, . 
\monend
We use the particular 
diffusive time step  $ \, \,  \Delta t = \Delta x ^2  $. 
Then 
$ \,   \sigma \equiv {{1}\over{s_J}}  -  {{1}\over{2}} = {1\over6} \, $ 
  and   the diffusivity  follows the relation 
  $  \, \kappa = {{\sigma_J}\over{3}} \, $ and 
\moneq \label{diffusivite}  
\kappa =  {1\over18} \, .  
\monend
%
%
We have chosen an odd number of mesh cells in these numerical experiments.
With the constraint $ \, \Delta t = \Delta x^2 $, it is not possible to 
obtain exactly the same exact final time. We have adapted the number of 
time   
steps in order to have very close values for the final time with the different meshes. 


\newpage 
\monitem 
Comparison with finite-difference approximation 

\noindent 
%
Remark that the solution of the  heat equation on a square with an initial Gaussian and 
periodic boundary conditions has 
to 
our knowledge no analytical solution. In consequence,
%
we compare the solution obtained by the lattice Boltzmann scheme with the result 
computed with two-dimensional finite differences, centered in space and explicit in time.
The degrees of freedom are located at half-integer  positions,   
exactly as done with the lattice Boltzmann scheme:
\moneqstar  
\rho_{ i+{1\over2} , \, j+{1\over2} }^{n} \approx \rho \Big( \Big( i+{1\over2} \Big) \, \Delta x , \, 
\Big( j+{1\over2} \Big) \, \Delta x , \,  n\, \Delta t \Big) \,. 
\monendstar 
We finite difference  the heat equation (\ref{chaleur}) in the following way~:
\moneqstar 
  \left \{ \begin {array}{l} 
\displaystyle {{1}\over{\Delta t}} \, \Big( 
\rho_{ i+{1\over2} , \, j+{1\over2}}^{n+1} - \rho_{i+{1\over2} , \, j+{1\over2} }^{n} \Big) 
 - \kappa \bigg[ 
{{1}\over{\Delta x^2}} \,\Big( 
\rho_{i+{3\over2} , \, j+{1\over2}}^{n} - 2 \, \rho_{ i+{1\over2} , \, j+{1\over2} }^{n} 
+ \rho_{i-{1\over2} , \, j+{1\over2} }^{n} \Big) \\ \displaystyle
\qquad  \qquad \qquad + 
{{1}\over{\Delta y^2}} \,\Big( 
\rho_{i+{1\over2} , \, j+{3\over2} }^{n} - 2 \, \rho_{i+{1\over2} , \, j+{1\over2} }^{n}
 + \rho_{i+{1\over2} , \, j-{1\over2} }^{n} \Big) \bigg]  = 0 \, . 
\end {array}  \right. \monendstar 
We use exactly the same grid in space for both schemes 
and exactly the same time step (and in consequence the same number of time steps). 
The parameters for both schemes are compared in Table~\ref{table-1}. 
%
\begin{table}  [H]     \centering
 \centerline { \begin{tabular}{|c|c|c|c|c|c|}    \hline 
 Number of cells & $ 13 \times 13 $  & $ 27 \times 27 $ & $ 55 \times 55 $  &  $ 111 \times 111 $   
& $ 223 \times 223  $  \\   \hline 
 nb. of time steps D2Q9  &  8 & 36 & 128 &  600 &  2048  \\   \hline  
 final time  &    $ 0.18935 $ & $ 0.19753 $ &  $ 0.17741 $ &  $ 0.19479 $ &   $ 0.16473  $ \\   \hline  
\end{tabular} }  
\caption{
D2Q9 numerical experiments with diffusive  scaling, $\, s_J \, $ given by (\ref{sJ-fixe}) 
and diffusivity $ \, \kappa \, $ by (\ref{diffusivite}). } 
\label{table-1} \end{table}
%

   \vspace {-.5cm}   \begin{figure}    [H]    \centering
 \includegraphics [ trim = {105mm 25mm 105mm 25mm}, clip, height=.35 \textwidth]
 {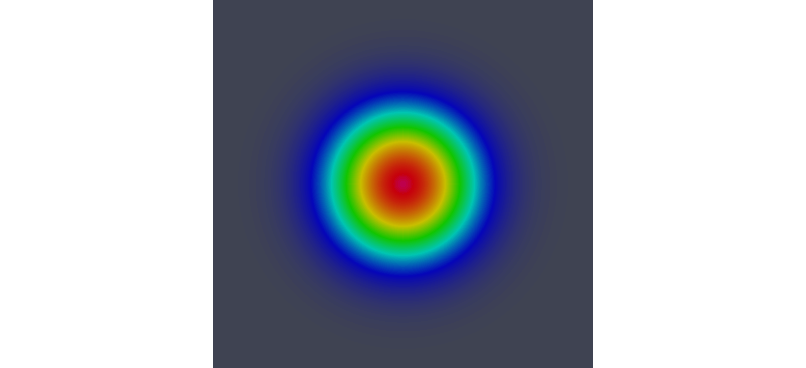}
\quad 
\includegraphics [ trim = {105mm 25mm 105mm 25mm}, clip, height=.35 \textwidth] 
{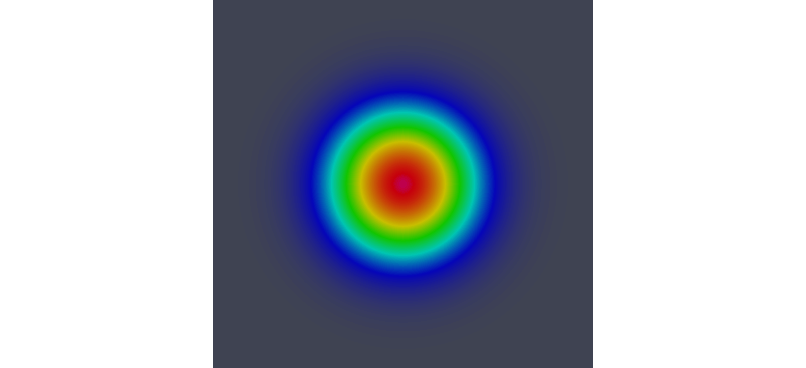}
\caption{Two-dimensional heat equation (\ref{chaleur}), 
$\, \kappa = {1\over18} $.  D2Q9 scheme with diffusive scaling (left) vs.  
explicit finite differences (right)~;  mesh 111 $\times$ 111, time = 0.19479. }
  \label{111-111-t019479} \end{figure}

   \vspace {-.5cm}   \begin{figure}    [H]  \centering
\includegraphics [ trim = {105mm 25mm 105mm 25mm}, clip, height=.35 \textwidth] 
{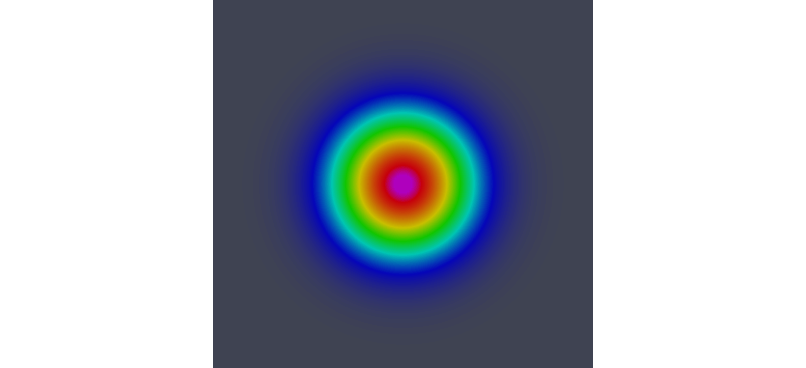}
\quad  
\includegraphics [ trim = {105mm 25mm 105mm 25mm}, clip, height=.35 \textwidth] 
{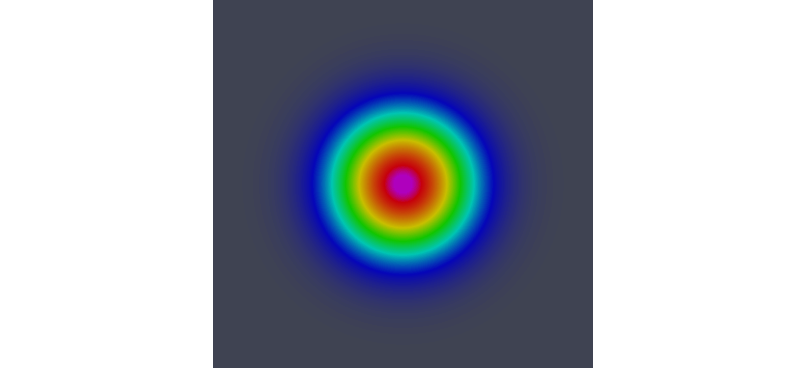}
\caption{Two-dimensional heat equation (\ref{chaleur}), 
$\, \kappa = {1\over18} $.  D2Q9 scheme with diffusive scaling (left) vs.
explicit finite differences (right)~;  mesh 223 $\times$ 223, time = 0.16473.  } 
  \label{223-223-t016473} \end{figure}

\begin{figure}    [H] \centering
\includegraphics [height=.55 \textwidth] {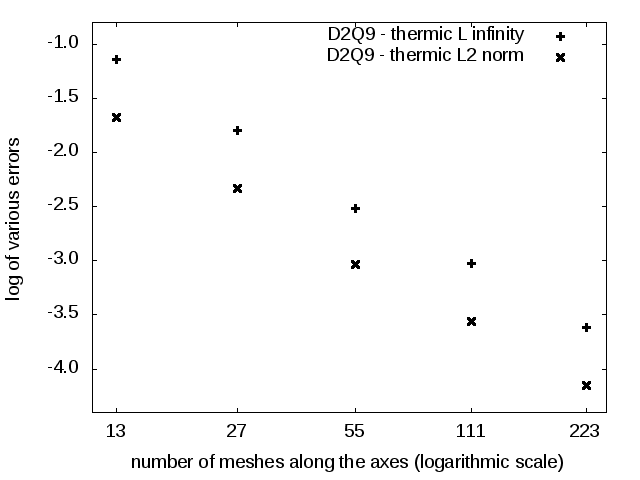}
\caption{Two-dimensional heat equation (\ref{chaleur}), 
$\, \kappa = {1\over18} $.
Difference of the numerical results computed with the 
 D2Q9 scheme with diffusive scaling and explicit finite differences
%
 at the final times presented in Table~1. 
%
The order of convergence for this residual in the  $ \, {\rm L}^\infty \, $ norm is equal to 3.41 
and  in the  $ \, {\rm L}^2 \, $ norm it is   3.96. }
  \label{convergence-diffusive-scaling} \end{figure}

\noindent 
The results follow what is expected. 
The approximate solutions 
 of both schemes are very similar as observed in Figures~\ref{111-111-t019479} 
and \ref{223-223-t016473} for $\, 111 \times 111 \, $ and $ \, 223 \times 223 \, $ meshes.
The difference between the two schemes 
%
at the final time 
%
is presented in  Figure~\ref{convergence-diffusive-scaling}. 
The order of convergence of   this residual is  approximately of order $4$. 
Since the finite difference method is  
of second order accuracy \cite{RM57}, this indicates that the lattice Boltzmann 
method approaches the heat equation with second-order  accuracy.   

   \vspace {-.5cm}   


\monitem 
Scalar D2Q9 numerical experiments with acoustic scaling

\noindent 
We still wish to solve the 
heat equation (\ref{chaleur}) 
in the square $ \, \Omega = \, [ -1 , \, 1 ]^2 $
 with periodic boundary conditions. 
The initial condition is again given by a Gaussian profile (\ref{gaussienne}). 
The given diffusivity  is imposed by the value (\ref{diffusivite}). 
We adopt an acoustic scaling with $\, \Delta t = \Delta x \, $ for the 
D2Q9 lattice Boltzmann simulations. 
We compare the results with explicit finite differences; 
in this case, we take $ \, \Delta t \simeq  \Delta x ^2 \, $  
and the time step is chosen in order to obtain exactly the same final time
than with the lattice Boltzmann method.

\begin{table}  [H]     \centering
 \centerline { \begin{tabular}{|c|c|c|c|c|c|}    \hline 
 Number of cells & $ 13 \times 13 $  & $ 27 \times 27 $ & $ 55 \times 55 $  &  $ 111 \times 111 $   
& $ 223 \times 223  $  \\   \hline 
  D2Q9 $\, s_J \,$ parameter &  1.5  & 1.182     & 0.830 &  0.52  &  0.298  \\   \hline 
 nb. of time steps D2Q9  &  8 & 16 & 32 &  64 &  128  \\   \hline  
 nb. of time steps, finite differences  &  8 & 32 & 128 &  512 &  2048  \\   \hline  
 final time  &   $ 0.18935 $ & $  0.18234 $ & $  0.17902  $ &  $ 0.17741 $ &  $ 0.17661  $  \\   \hline  
\end{tabular} }  
\caption{
D2Q9 numerical experiments with acoustic scaling. The diffusivity $ \, \kappa = {1\over18} \,$
is imposed in all the simulations.  } \label{table-2} \end{table}

   \vspace {-.5cm}   \begin{figure}    [H]    \centering
 \includegraphics [ trim = {105mm 25mm 105mm 25mm}, clip, height=.35 \textwidth] 
{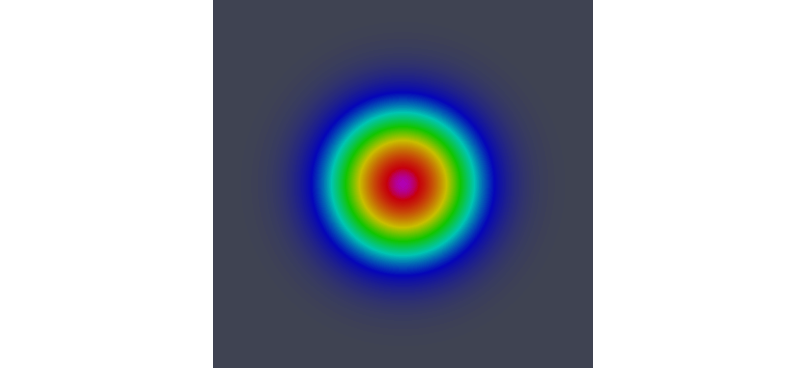}
\quad 
\includegraphics [ trim = {105mm 25mm 105mm 25mm}, clip, height=.35 \textwidth]
 {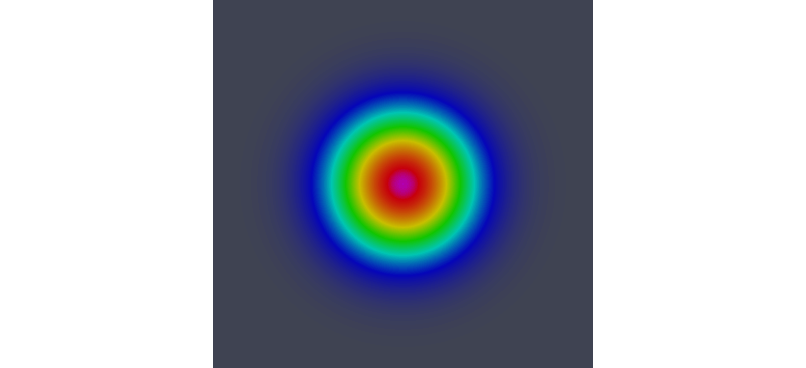}
\caption{Two-dimensional heat equation (\ref{chaleur}), 
$\, \kappa = {1\over18} $. D2Q9   lattice Boltzmann scheme  with acoustic scaling (left) vs. 
explicit finite differences (right)~;  results at time =  0.17741  for a 
111 $\times$ 111 mesh. }
  \label{111-111-acoustic-t017741} \end{figure}

   \vspace {-.5cm}   \begin{figure}    [H] \centering
 \includegraphics [ trim = {105mm 25mm 105mm 25mm}, clip, height=.35 \textwidth] 
 {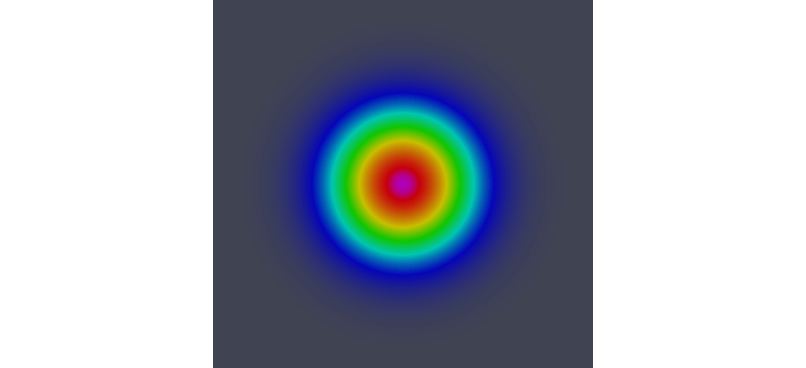}
\quad 
 \includegraphics [ trim = {105mm 25mm 105mm 25mm}, clip, height=.35 \textwidth]
 {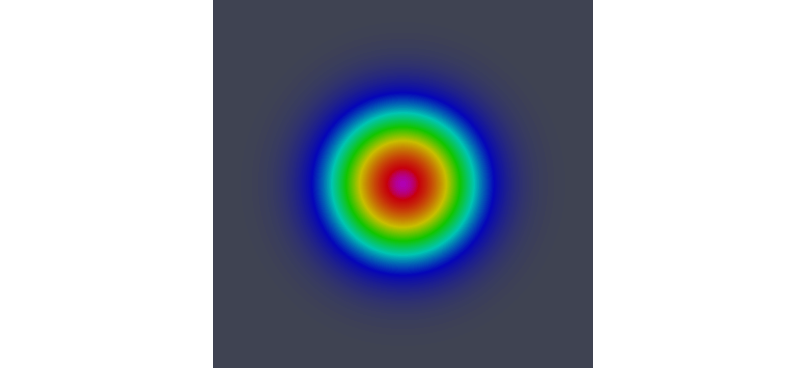}
\caption{Two-dimensional heat equation (\ref{chaleur}), 
$\, \kappa = {1\over18} $. D2Q9 lattice Boltzmann scheme with acoustic scaling (left) vs. 
explicit finite differences (right)~;  results at time = 0.1766 for a 
223 $\times$ 223 mesh. }
  \label{223-223-acoustic-t01766} \end{figure}

\begin{figure}    [H] \centering
\includegraphics [height=.55 \textwidth] {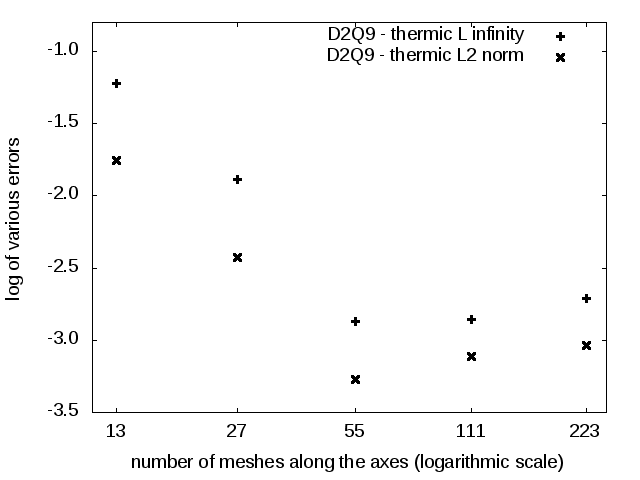}
\caption{Two-dimensional heat equation (\ref{chaleur}), 
$\, \kappa = {1\over18} $, D2q9 scheme with acoustic scaling  vs. 
explicit finite differences 
%
 at the final times presented in Table~2. 
%
There is no numerical evidence of coherence between the two methods     
when the mesh is refined. }
  \label{no-convergence-acoustic-scaling} \end{figure}

\monitem
The numerical results presented in Figures \ref{111-111-acoustic-t017741}
and \ref{223-223-acoustic-t01766}  for the two meshes of 111 $\times $ 111 and\br
  223 $\times $ 223 seem  correct. 
But a quantitative examination of the results (Figure~\ref{no-convergence-acoustic-scaling})  
shows that after a convergence similar to the one obtained for diffusive scaling
(see Figure~\ref{convergence-diffusive-scaling}), a  persistent difference
appears. This qualitative behaviour is very similar to what has been 
observed in \cite{BDGLT18}   in one spatial  dimension.

   \vspace {-.5cm}   
%

\newpage 
A new  analysis of the scheme is necessary to explain this lack of convergence towards the 
expected diffusive model.


\bigskip \bigskip   \noindent {\bf \large    4) \quad 
Dispersion  equation for an evanescent  relaxation  } 

\noindent 
In this section, we propose a first-order  analysis when the last relaxation coefficients in 
(\ref{coef-relax-juillet2017}) remain 
 fixed or  when the relaxation coefficient $ \, s_J \, $ 
for the momentum $ \, J \, $ follows the choice presented in Eq.~(\ref{s-kappa}), {\it id est}
\moneq \label{sJ-delta-t} 
s_J = {{4+\alpha}\over{6}} \, \lambda^2 \, {{\Delta t}\over{\kappa}} + {\rm O} (\Delta t^2) \, . 
\monend 
%

\monitem Fixed  relaxations

\noindent 
We write the relation (\ref{LB-schema}) in terms of the moments $ \, m \, $ defined in 
(\ref{moments}): 
\moneq \label{iter-moments}  
m_k(x, \, t + \Delta t) = \sum_{j \, \ell} M_{k \, j} \, M^{-1}_{j \, \ell} \, m_\ell^*(x - v_j \, \Delta t , \, t ) 
\, . 
\monend
Before doing a Taylor expansion at order 1, we introduce the following ``momentum velocity'' operator
matrix $ \, \Lambda \, $ defined according to 
\moneq \label{grand-lambda}  
\Lambda_{k \, \ell}  \equiv -  \sum_{j \, \ell \, \alpha} M_{k \, j} \, v_j^\alpha \, M^{-1}_{j \, \ell}  
\, {{\partial}\over{\partial x_\alpha}}  \, . 
\monend
For the D2Q9 scheme, this matrix can be explicitly calculated  \cite{Du09} and we have  
\moneq \label{grand-lambda-d2q9}  
\Lambda = - \left(\begin{array}{|c|cccccccc|}
\hline  
0 & \partial_x & \partial_y & 0 & 0 & 0 & 0 & 0 & 0 \\ \hline  
{{2\lambda^2}\over{3}} \, \partial_x & 0 & 0 & {1\over6} \, \partial_x &  {1\over2} \, \partial_x &
 \partial_y & 0 & 0 & 0 \\
{{2\lambda^2}\over{3}} \, \partial_y & 0 & 0 & {1\over6} \, \partial_y &  -{1\over2} \, \partial_y &
 \partial_x & 0 & 0 & 0 \\
0 & \lambda^2 \, \partial_x & \lambda^2 \, \partial_y & 0 & 0 & 0 &  \partial_x &  \partial_y & 0 \\ 
0 &  {{\lambda^2}\over{3}} \, \partial_x &  -{{\lambda^2}\over{3}} \, \partial_y & 0 & 0 & 0 &  
- {1\over3} \, \partial_x &  {1\over3} \, \partial_y & 0 \\ 
0 &  {{\lambda^2}\over{3}} \, \partial_y &  -{{\lambda^2}\over{3}} \, \partial_x & 0 & 0 & 0 &  
- {1\over3} \, \partial_y &  {1\over3} \, \partial_x & 0 \\ 
0 & 0 & 0 &   {{\lambda^2}\over{3}} \, \partial_x & - \lambda^2 \, \partial_x & \lambda^2 \, \partial_y &
0 & 0 & {1\over3} \, \partial_x  \\ 
0 & 0 & 0 &   {{\lambda^2}\over{3}} \, \partial_y &  \lambda^2 \, \partial_y & \lambda^2 \, \partial_x &
0 & 0 & {1\over3} \, \partial_y \\ 
0 & 0 & 0 & 0 & 0 & 0 & \lambda^2 \, \partial_x & \lambda^2 \, \partial_y & 0 \\ \hline  
\end{array}\right)  \, . \monend
%
We split the moment vector into two blocks:   
\moneq \label{moments-2-composantes}  
m =  \begin{pmatrix} W \\ Y 
\end{pmatrix}  \monend
with $ \, W = \rho \, $ in our scalar  example and $ \, Y \, $ a column vector with 8 components. 
%
We decompose also the operator matrix $ \, \Lambda \, $ into four blocks that respect   
the decomposition (\ref{moments-2-composantes}): 
\moneq \label{grand-lambda-abcd}  
\Lambda \equiv \begin{pmatrix} A & B \\ C & D \end{pmatrix} \, . 
\monend    
In our case, $ \, A \, $ is a scalar $ \, 1 \times 1 \, $ matrix, 
$ \, B \,$ has one line and 8 columns, 
$ \, C \, $  is composed by 8 lines and 1 column and  $ \, D \, $ is a  8$\times$8 square matrix
as shown in the right-hand  side of relation~(\ref{grand-lambda-d2q9}). 
We can also introduce a constant matrix $ \, E \, $ with 8  lines and one column such that 
the relation (\ref{moments-equilibre}) can be written in the form 
\moneq    \label{def-E} 
Y^{\rm eq}  \, \equiv  \, E \, m  \, . 
\monend    
The relation (\ref{iter-moments}) is expanded at first order:
\moneq \label{lbm-ordre-1}  
m + \Delta t \, \partial_t m + {\rm O}(\Delta t^2) = m^* + \Delta t \,  \Lambda \, m^* + {\rm O}(\Delta t^2) 
 \monend
and due to (\ref{def-E}), we have 
\moneq \label{m-star}  
 m^* = \begin{pmatrix} {\rm I} & 0 \\ S \, E &  {\rm I}-S \end{pmatrix}  \, m \, . \monend
%
The relation (\ref{lbm-ordre-1}) can be written in the form
\moneq \label{L-m-ordre-2}  
L \, m \equiv m^* - m + \Delta t \, \, \big( - \partial_t m + \Lambda \,  m^* \, \big)   
\, = \,   {\rm O}(\Delta t^2) \,, 
 \monend
with 
\moneq \label{matrice-L}  
L \equiv \begin{pmatrix} 0 & 0 \\ S \, E & -S \end{pmatrix} + 
\Delta t \, \left[ \begin{pmatrix} -\partial_t & 0 \\ 0 & -\partial_t  \end{pmatrix}
+ \begin{pmatrix} A & B \\ C & D \end{pmatrix} \,\, 
\begin{pmatrix}  {\rm I} & 0 \\ S \, E &  {\rm I}-S \end{pmatrix}  \right] \, .  
 \monend
%

\noindent  The dispersion relation associated with the relation (\ref{L-m-ordre-2})  
can be written in a simple way:
\moneq \label{dispersion-L}  
{\rm det} \, L = 0 \, . 
 \monend
We expand this determinant in order to eliminate the non-conserved moments $\, Y $. 
Moreover, due to the right-hand  side of Eq.~(\ref{L-m-ordre-2}), 
we can neglect all the terms of second or third order relative to $ \, \Delta t $. 
We write the expression (\ref{matrice-L}) of the matrix $ \, L \, $ in  the form 
\moneqstar 
L =  \begin{pmatrix} \Delta t \, (-\partial_t + A + B \, S \, E ) & \Delta t \, B \, (  {\rm I}-S ) \\
S \, E + \Delta t \, ( C + D \, S \, E ) & -S +  \Delta t \, ( -\partial_t + D \, (  {\rm I}-S ) ) 
\end{pmatrix} \, . \monendstar 
We apply Gaussian elimination in order to make explicit the condition (\ref{dispersion-L}). 
We multiply this matrix at left by the regular  matrix $ \, K \, $ defined by 
\moneq \label{matrice-K}  
K =  \begin{pmatrix}  {\rm I} & \Delta t \, B \, (  {\rm I}-S ) \, S^{-1} \\ 0 &  {\rm I} 
\end{pmatrix} \, .  \monend
Then we have, after some lines of algebra, 

\smallskip \noindent
$ \displaystyle K \, L =  
\begin{pmatrix}  {\rm I} & \Delta t \, B \, (  {\rm I}-S ) \, S^{-1} \\ 0 &  {\rm I} \end{pmatrix} 
\,\,  \begin{pmatrix} \Delta t \, (-\partial_t + A + B \, S \, E ) & \Delta t \, B \, (  {\rm I}-S ) \\
S \, E + \Delta t \, ( C + D \, S \, E ) & -S +  \Delta t \, ( -\partial_t + D \, (  {\rm I}-S ) ) 
\end{pmatrix} $ 

\smallskip \noindent $ \qquad =  
\begin{pmatrix} \Delta t \, (-\partial_t + A + B \, S \, E ) +  
 \Delta t \, B \, (  {\rm I}-S ) \, S^{-1} \, S \, E & {\rm O} (\Delta t^2 ) \\ 
S \, E + {\rm O} (\Delta t ) & -S + {\rm O} (\Delta t ) 
\end{pmatrix} $ 

\smallskip \noindent
and we have the following triangular form for the product $ \, K \, L $:
\moneqstar 
 K \, L = \begin{pmatrix}  \Delta t \, (-\partial_t + A + B \, E ) &  {\rm O}(\Delta t^2) \\
S \, E +  {\rm O}(\Delta t) & -S +  {\rm O}(\Delta t) 
\end{pmatrix} \, . \monendstar 
%
Then the relation (\ref{dispersion-L}) is equivalent at first order to the following set   
of first order partial differential equations:
\moneq \label{edp-ordre-1}  
( -\partial_t + A + B \, E ) \, W =  {\rm O}(\Delta t) \, , \monend
recovering the first step of the Berlin algorithm presented in Augier {\it et al.} \cite {ADG14}.  
For the scalar diffusion problem, this equation expresses  simply that 
\moneqstar 
 \partial_t \rho =  {\rm O}(\Delta t) \, . 
\monendstar 
This result is 
consistent   
 with the second-order  analysis presented at the relation  in (\ref{edp-diff}).
%

When we use diffusive scaling, this dispersion  equation can be adapted in order 
to recover the heat equation at zero order of accuracy. It is then equivalent to 
the Taylor expansion method with the diffusive scaling, as used in  \cite{DL13}.  

\newpage 
\monitem Evanescent  relaxations

\noindent 
When $ \, \Delta t \, $ and $ \, \Delta x \, $ tend to zero with the acoustic scaling, 
these  two infinitesimals are of the  same order. 
The expansion (\ref{sJ-delta-t}) of the relaxation coefficient $ \, s_J \, $ implies 
that the previous asymptotic calculus has to be made more precise. 
The coefficient $ \, s_J \, $ is now at first order proportional to the time step $ \, \Delta t $. 
We decompose the non-conserved  moments $ \, Y \, $ into two families:
the quasi-conserved  moments $ \, U \,$ {\it id est} the two
components of the momentum $ \, J \, $ in the scalar  case-- and the other truly
non-conserved  moments $\, Z $:
\moneq \label{decomposition-Y}  
Y =  \begin{pmatrix} U \\ Z 
\end{pmatrix}  \, . \monend
%
%
The 8-component vector $ \, Y \, $ is split into a first vector $ \, U \in \RR^2 \, $   
and a second one $ \, Z \, $ with 6~components. 
In other words, the family of moments is split into three components:   
\moneqstar
m =  \begin{pmatrix} W \\ U \\ Z 
\end{pmatrix}  \, . \monendstar
Then the 8$\times$8 relaxation matrix $ \, S \, $ can be decomposed into two blocks:   
\moneq \label{decomposition-S}  
S =  \begin{pmatrix} \Delta t \, {\widetilde S} + {\rm O}(\Delta t^2) & 0  \\ 
0 & S_Z \end{pmatrix}  \, . \monend
The top left block in the right hand side of   
 (\ref{decomposition-S}) tends to zero as the mesh is refined. 
The equilibrium vector $ \, E \, $ is naturally split into the  quasi-conserved    
component $ \, E_U \, $ and the truly relaxing component $ \, E_Z $:       
\moneq \label{decomposition-E}  
E =  \begin{pmatrix} E_U \\ E_Z \end{pmatrix}  \, . \monend
We have:  
$\,\, \displaystyle S \, E = 
\begin{pmatrix} \Delta t \, {\widetilde S} & 0  \\ 0 & S_Z \end{pmatrix} \,\, 
\begin{pmatrix} E_U \\ E_Z \end{pmatrix} 
\, = \, \begin{pmatrix}  \Delta t \, {\widetilde S} \, E_U \\  S_Z \,  E_Z \end{pmatrix} \,\, $ and 
the relation (\ref{m-star})  takes the form 

\moneq \label{decomposition-m-star}  
m^*  =  \begin{pmatrix}  {\rm I} & 0 & 0 \\ 
\Delta t \, {\widetilde S} \, E_U &   {\rm I} - \Delta t \, {\widetilde S} & 0 \\
 S_Z \,  E_Z & 0 &  {\rm I} - S_Z  \end{pmatrix} \,\, m  \, . \monend
%
Then the momentum velocity operator matrix $ \, \Lambda \, $ is split into 9 blocks:  
\moneq \label{decomposition-grand-lambda}  
 \Lambda =  \begin{pmatrix} A & A_2 & B_1 \\ A_3 & A_4 & B_2 \\ C_1 & C_2 & D_4  \end{pmatrix}  \, .\monend
This block  structure (\ref{decomposition-grand-lambda}) is explicitly given for our thermal D2Q9 
in  the form 
\moneqstar
\Lambda = - \left(\begin{array}{|c|cc|cccccc|} \hline 
 0 & \partial_x & \partial_y & 0 & 0 & 0 & 0 & 0 & 0 \\ \hline 
{{2\lambda^2}\over{3}} \, \partial_x & 0 & 0 & {1\over6} \, \partial_x &  {1\over2} \, \partial_x &
 \partial_y & 0 & 0 & 0 \\
{{2\lambda^2}\over{3}} \, \partial_y & 0 & 0 & {1\over6} \, \partial_y &  -{1\over2} \, \partial_y &
 \partial_x & 0 & 0 & 0 \\ \hline 
0 & \lambda^2 \, \partial_x & \lambda^2 \, \partial_y & 0 & 0 & 0 &  \partial_x &  \partial_y & 0 \\ 
0 &  {{\lambda^2}\over{3}} \, \partial_x &  -{{\lambda^2}\over{3}} \, \partial_y & 0 & 0 & 0 &  
- {1\over3} \, \partial_x &  {1\over3} \, \partial_y & 0 \\ 
0 &  {{\lambda^2}\over{3}} \, \partial_y &  -{{\lambda^2}\over{3}} \, \partial_x & 0 & 0 & 0 &  
- {1\over3} \, \partial_y &  {1\over3} \, \partial_x & 0 \\ 
0 & 0 & 0 &   {{\lambda^2}\over{3}} \, \partial_x & - \lambda^2 \, \partial_x & \lambda^2 \, \partial_y &
0 & 0 & {1\over3} \, \partial_x  \\ 
0 & 0 & 0 &   {{\lambda^2}\over{3}} \, \partial_y &  \lambda^2 \, \partial_y & \lambda^2 \, \partial_x &
0 & 0 & {1\over3} \, \partial_y \\ 
0 & 0 & 0 & 0 & 0 & 0 & \lambda^2 \, \partial_x & \lambda^2 \, \partial_y & 0  \\ \hline 
\end{array}\right)  \, . \monendstar
Then 

\smallskip \noindent
$  \displaystyle L = \begin{pmatrix} 0 & 0 & 0 \\ 
\Delta t \, {\widetilde S} \, E_U &  - \Delta t \, {\widetilde S} & 0 \\
 S_Z \,  E_Z & 0 &   - S_Z  \end{pmatrix} 
- \Delta t  \, \begin{pmatrix} \partial_t  & 0 & 0 \\ 0 & \partial_t  & 0 \\ 0 & 0 & \partial_t \end{pmatrix} $

\qquad \qquad  \qquad  
$ \displaystyle  +   \Delta t  \,
\begin{pmatrix} A & A_2 & B_1 \\ A_3 & A_4 & B_2 \\ C_1 & C_2 & D_4  \end{pmatrix} \,\, 
 \begin{pmatrix}  {\rm I} & 0 & 0 \\ 
\Delta t \, {\widetilde S} \, E_U &   {\rm I} - \Delta t \, {\widetilde S} & 0 \\
 S_Z \,  E_Z & 0 &  {\rm I} - S_Z  \end{pmatrix} \, . $

\smallskip \noindent  
This expression can be expanded to  first order in  $ \, \Delta t \, $
without any change in the result of the Gaussian elimination. Then 
we can neglect the terms of order one in $ \, \Delta t \, $ in the 
last product of two matrices. We obtain 

\smallskip \noindent 
$  \displaystyle 
\begin{pmatrix} A & A_2 & B_1 \\ A_3 & A_4 & B_2 \\ C_1 & C_2 & D_4  \end{pmatrix} \,\,  
 \begin{pmatrix}  {\rm I} & 0 & 0 \\ 0  &   {\rm I}  & 0 \\
 S_Z \,  E_Z & 0 &  {\rm I} - S_Z  \end{pmatrix}  \,=\, 
\begin{pmatrix} A + B_1 \,  S_Z \,  E_Z & A_2 & B_1 \, ( {\rm I} - S_Z ) \\
A_3 + B_2 \,  S_Z \,  E_Z & A_4 & B_2 \,  ( {\rm I} - S_Z ) \\
C_1 + D_4 \,  S_Z \,  E_Z & C_2 &  D_4 \,  ( {\rm I} - S_Z )  \end{pmatrix}   \,, $ 

\smallskip \noindent 
and, up to order $ \, {\rm O }(\Delta t) $, we have 
\moneq \label{decomposition-L}  
L \! = \! \begin{pmatrix} \Delta t  ( - \partial_t + A + B_1 \, S_Z \, E_Z ) & 
\!\!\! \Delta t \, A_2 & 
\!\!\!  \Delta t \, B_1 \, ( {\rm I} - S_Z )  \\
 \Delta t  (  {\widetilde S} \, E_U + A_3 + B_2 \,  S_Z \,  E_Z ) & 
\!\!\! \Delta t  ( - {\widetilde S}  - \partial_t  + A_4 ) & 
\!\!\!   \Delta t \, B_2 \,  ( {\rm I} - S_Z ) \\
 S_Z \, E_Z + \Delta t \, ( C_1 + D_4 \,  S_Z \,  E_Z ) & 
\!\!\! \Delta t \, C_2 &
\!\!\!\!\!\!\!\! -S_Z + \Delta t  ( - \partial_t + D_4  ( {\rm I} - S_Z )) 
\end{pmatrix} \! . \monend

\noindent 
With the method of Gaussian elimination used previously, 
we  multiply the matrix $ \, L \, $ obtained in (\ref{decomposition-L}) 
on the left by the following matrix
\moneqstar K' = \begin{pmatrix} {\rm I} & 0 & \Delta t \, B_1 \, ( {\rm I} - S_Z ) \, S_Z^{-1} \\
0 & {\rm I} & \Delta t \, B_2 \, ( {\rm I} - S_Z ) \, S_Z^{-1} \\
0 & 0 &  {\rm I} \end{pmatrix} \monendstar 
whose determinant is equal to 1. 
After some  elementary algebra, we obtain 
\moneqstar K' \, L = 
\begin{pmatrix}  \Delta t \,  (- \partial_t + A + B_1 \,  E_Z ) & \Delta t \, A_2 & {\rm O}(\Delta t^2) \\
 \Delta t \,  (  {\widetilde S} \, E_U + A_3 + B_2 \,  E_Z ) & 
\Delta t \,  (- {\widetilde S}  - \partial_t  + A_4 ) & {\rm O}(\Delta t^2) \\
S_Z \, E_Z +  {\rm O}(\Delta t) &  \Delta t \, C_2 & - S_Z +  {\rm O}(\Delta t) 
 \end{pmatrix} \, . \monendstar 
On one hand,  $ \, {\rm det} K' = 1 \, $ and on the other hand, 
the last column of the matrix $ \, K' \, L \, $ 
is composed of negligible terms except for the last one.   
Then  we have the condition  (\ref{dispersion-L}) if and only if 
the determinant of the $2\times 2$ upper block matrix is null. In other terms,   
this matrix has a nontrivial kernel at order~one relative  to $\, \Delta t \, $ and 
we have 
\moneq \label{edp-matriciel}  
\begin{pmatrix}  - \partial_t + A + B_1 \,  E_Z  & A_2 \\
  {\widetilde S} \, E_U + A_3 + B_2 \,  E_Z &  - \partial_t - {\widetilde S}   + A_4 \end{pmatrix} \, 
\begin{pmatrix} W \\ U \end{pmatrix} =  {\rm O}(\Delta t)  \, . 
\monend
Then the equivalent partial differential equations are written as a system
involving the conserved variable $ \, W \, $ and the quasi conserved moments $\, U $:
\moneq \label{system-edp}  
  \left \{ \begin {array}{l}
\partial_t W = (  A + B_1 \,  E_Z ) \, W + A_2 \, U +  {\rm O}(\Delta t)  
\\  \vspace{-4 mm} \\ \displaystyle 
\partial_t U +  {\widetilde S} \, U = ( A_3 + B_2 \,  E_Z +   {\widetilde S} \, E_U ) \, W 
+ A_4 \, U +  {\rm O}(\Delta t) \, . 
\end{array} \right. \monend
This result generalizes the first analysis done in \cite{BDGLT18}   
for the D1Q3 scheme. 
%
When we replace the block matrices  introduced in the relations   
(\ref{decomposition-S}), (\ref{decomposition-E}) and 
(\ref{decomposition-grand-lambda}) 
by their D2Q9 values, we establish   that with the acoustic scaling, 
the  scalar  D2Q9 lattice Boltzmann scheme with acoustic scaling 
admits the following asymptotic damped acoustic model 
\moneq \label{damped-acoustic} 
  \left \{ \begin {array}{l}
\displaystyle  {{\partial \rho}\over{\partial t}}  + {\rm div} J = {\rm O}(\Delta x) 
\\  \vspace{-4 mm} \\\displaystyle
 {{\partial J_\alpha }\over{\partial t}}   +  c_0^2 \,\, 
{{\partial \rho}\over{\partial x_\alpha}} + g  \, J_\alpha = {\rm O} (\Delta x) \, , \,\, 
1 \leq \alpha \leq 2 \,, 
\end{array} \right. \monend
with a sound velocity $ \, c_0 \, $ and a damping coefficient $ \, g \, $ 
given by the relations 
\moneq \label{c0-gg}  
 c_0^2 = {{\lambda^2}\over{6}} \, (4 + \alpha) \,, \quad 
g =  {{c_0^2}\over{\kappa}}  \, . 
\monend
The above is a very interesting analysis, 
and clearly the correct two-dimensional analog of the earlier result for D1Q3.  
We point out that it is equivalent to a damped wave equation. 

\bigskip \bigskip   \noindent {\bf \large    5) \quad  
Scalar  D2Q9 scheme converging towards  damped acoustic  }   

\noindent 
We have now two partial differential equations with which  to compare the numerical solution obtained
with the scalar  D2Q9  lattice Boltzmann scheme: the initial heat equation (\ref{chaleur})
and the damped acoustic system (\ref{damped-acoustic}). We first consider numerical experiments
done in Section~3 and compare our previous results with this new model. 
We also study  in  detail the eigenmodes of the system  (\ref{damped-acoustic}) and propose
a simple numerical experiment with a  sinusoidal analytic solution. 
The evolution of an initial Gaussian is again performed, with two diffusion 
coefficients varying by one order of magnitude. 

\monitem 
Damped acoustics as a limiting  model for the previous numerical experiments?   

\noindent 
We  wish to approximate the  system of 
damped acoustic  equations (\ref{damped-acoustic}). 
The sound velocity is given by  (\ref{c0-gg}). With the choice 
(\ref{alpha-beta}), we obtain the classical value $ \, c_0 = {{\lambda}\over{\sqrt{3}}} $. 
The imposed diffusivity $ \, \kappa \, $ and the relation  (\ref{c0-gg}) fix 
the value $ \, g = 6 \, $ 
for the zero-order  damping in the momentum equation of  (\ref{damped-acoustic}).
The geometry is the square $ \, \Omega = \, [ -1 , \, 1 ]^2 $
 with periodic boundary conditions. 
The initial density is still given by a Gaussian profile (\ref{gaussienne}).
Because the momentum $ \, J \, $ at equilibrium is identically null, 
we have taken this specific value as 
initial condition of our lattice Boltzmann simulations. We
suppose in consequence that the  initial condition for the momentum is 
simply $ \, J(x, t=0) = 0 $. 

\noindent 
We adopt acoustic scaling with $\, \Delta t = \Delta x \, $ for the 
D2Q9 lattice Boltzmann simulations. 
For the acoustic system (\ref{damped-acoustic}), 
we use explicit finite differences with staggered grids, hereafter named as ``HaWAY''  method
and described with some details in the Appendix. 
In this case, the acoustic time step  $ \, \Delta t_a \, $ for finite-difference   
simulations is proportional to the spatial  step  $ \, \Delta x $, 
with a stability constraint. 
The corresponding experiments are described in Table~\ref{table-3}.

\begin{table}  [H]     \centering
 \centerline { \begin{tabular}{|c|c|c|c|c|c|}    \hline 
 number of cells &  13 $\times$ 13  &  27 $\times$ 27  &  55 $\times$ 55   
&   111 $\times$ 111  &  223 $\times$ 223  \\   \hline 
D2Q9  $\, s_J \,$ parameter  &  1.5 &  1.182 &  0.830 &   0.52 &  0.298    \\   \hline 
 nb. of time steps d2q9 &   8 & 16 &  32 & 64 &  128   \\   \hline  
 idem, finite differences  &  32 & 64 &  128 & 256 &  512  \\   \hline  
final time   &   0.18935  & 0.18234 &  0.17902 &  0.17741  & 0.17661  \\   \hline  
\end{tabular} }  
\caption { Numerical experiments with the D2Q9 lattice Boltzmann test case
studied in Section~3 compared with the damped acoustic model (\ref{damped-acoustic})
simulated with the HaWAY   method.   } \label{table-3} \end{table}

   \vspace {-.5cm}   \begin{figure}    [H] \centering
\includegraphics [height=.55 \textwidth] {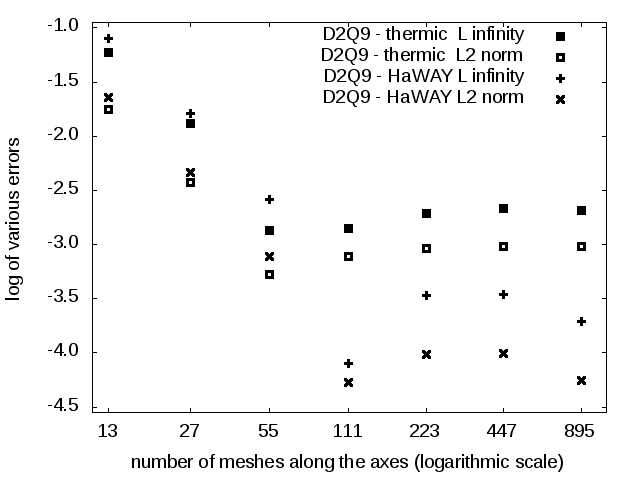}
\caption{Numerical  results for the damped acoustic model
with the experimental plan proposed in Table ~\ref{table-3}. 
}   
  \label{no-convergence-acoustic-wave} \end{figure}

\noindent 
The results obtained with this new experiment are very similar to the one 
obtained in Section~3. In particular, the numerical results computed with the 
damped acoustic model are very close to the ones presented   
in  Figures \ref{111-111-acoustic-t017741}
and \ref{223-223-acoustic-t01766}. 
When we look to the convergence with quite fine grids (Figure~\ref{no-convergence-acoustic-wave}), 
the signal is better than in  Figure~\ref{no-convergence-acoustic-scaling} but    
this experiment is still not entirely convincing.    

\monitem 
 Waves for the damped acoustic model 

\noindent 
We search modes of the type 
\moneq \label{modes-damped-acoustic} 
  \left \{ \begin {array}{l}
\displaystyle \rho = \rho_0 \, {\rm exp} \, ( -\gamma \, t + i \, k \cdot  x )
\\  \vspace{-4 mm} \\ \displaystyle 
J  = J_0 \, {\rm exp} \, (  -\gamma \, t +  i \, k \cdot x )   
\end{array} \right. \monend
for the damped acoustic model (\ref{damped-acoustic}--\ref{c0-gg}). Then 
we have to solve the following ill-posed linear system:
\moneq \label{3-modes} 
\displaystyle  - \gamma \,  \rho_0 + i \, k \cdot J_0 = 0 \,, \quad 
 i \, k \cdot \rho_0 + ( g - \gamma) \, J_0 = 0 \, .  
\monend
A first solution is  a transverse stationary wave 
with  $ \, \gamma = g $, $ \,\, \rho_0 = 0 \, $ and $ \,  k \cdot J_0 = 0 $. 
We do not consider this mode in this contribution. 
Then the other modes satisfy the following dispersion relation 
\moneq \label{dispersion-acoustic} 
\gamma^2 - g \, \gamma + | k |^2 \, c_0^2 = 0 \, . 
\monend
This equation has complex propagative roots when  
\moneq \label{propagative-acoustic} 
 g < 2 \,  | k | \, c_0  \,,
\monend
i.e.,  when the diffusivity $ \, \kappa \, $ is suficiently large
measured in a scale system based on the sound velocity and wave number:   
\moneqstar 
\kappa > {{c_0}\over{2 \,  | k | }}  \, . 
\monendstar 
In that case, the eigenvalue $ \, \gamma \, $ takes the form
\moneq \label{gamma-acoustic} 
 \gamma = {{g}\over{2}} \mp i \, \omega \,, \quad 
\omega = \sqrt { | k |^2 \, c_0^2 - {{g^2}\over{4}} } \, . 
\monend
The eigenvectors are finally given according to 
\moneq \label{vecteurs-propres}  \left \{ \begin {array}{l} 
\displaystyle  \rho = \rho_0 \, {\rm exp} \, \Big( -{{g}\over{2}} \, t \Big) \,\,
 {\rm exp} \, \big( i (  k \cdot  x \pm \, \omega \, t ) \big)  
  \\  \vspace{-4 mm} \\ \displaystyle 
J = i \, {{k}\over{| k |^2}} \, \rho_0 \, \Big(  -{{g}\over{2}}  \mp  i \, \omega \Big) \, 
{\rm exp} \, \Big( -{{g}\over{2}} \, t \Big) \,\, 
 {\rm exp} \, \big( i (  k \cdot  x \pm  \, \omega \, t ) \big) \, .  
  \end{array} \right.  \monend
%
We consider a pure analytical test case as the next experiment. 

\monitem 
 A two-dimensional sinusoidal wave 
 
\noindent
We keep  the value  $ \, \kappa = {1\over18}  \simeq 0.05555 \, $  
of the diffusivity introduced in (\ref{diffusivite}). 
We use the traditional value $ \, c_0 = {{1}\over\sqrt{3}} \, $ 
and the dissipation coefficient $ \, g \, $ (see (\ref{c0-gg})) is still 
equal to $ \, g = 6 $.   
We change the domain and consider $ \, [ 0 ,\, 2 \, \pi ]  ^2 \, $ with the initial condition 
$ \, \rho = \cos \big( ( 2  \pi \, (x+y) \big) \, $ and $ \, J = 0 $. 
Then $ \, k = 2 \, \pi \, (1, \, 1 ) \,  $ and the right-hand  side of 
(\ref{propagative-acoustic}) is $ \,  2 \,  | k | \, c_0 \simeq 5.924 $. 
In this case, the damped acoustic model  (\ref{damped-acoustic}) 
 exhibits a non-propagative mode.

\begin{figure}    [H] \centering
\includegraphics [trim = {100mm 50mm 100mm 50mm}, clip, height=.44 \textwidth] 
{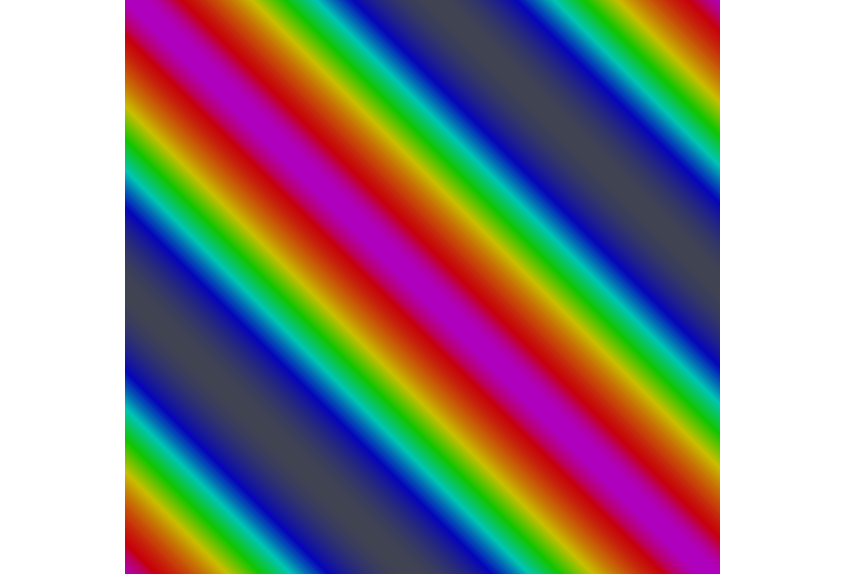}
\caption{  Two-dimensional wave with wave vector $ \, k = 2  \pi \, (1, \, 1 ) $. Initial 
condition.  }
  \label{d2q9periodic-densite-m021280} \end{figure}

\begin{figure}    [H] \centering
\includegraphics [height=.55 \textwidth] {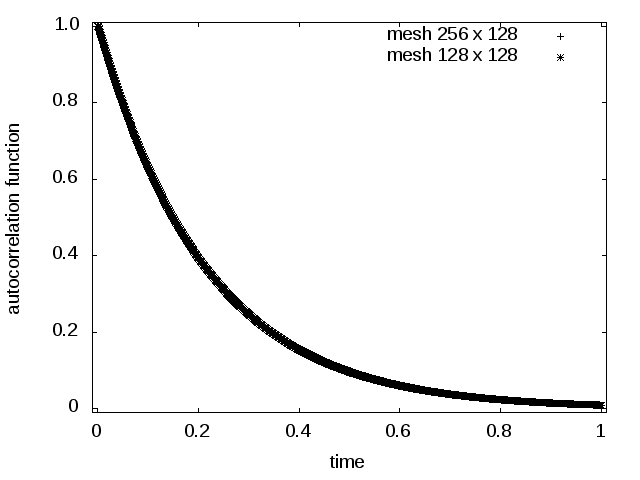}
\caption{  Two-dimensional wave with wave vector $ \, k = 2  \pi \, (1, \, 1 ) $, 
$ \, g \simeq  6  $, $ \,  2 \,  | k | \, c_0 \simeq 5.924 $. 
Autocorrelation of density.  }
  \label{manip02long-autocorrelation} \end{figure}

\begin{figure}    [H] \centering
\includegraphics [height=.55 \textwidth] {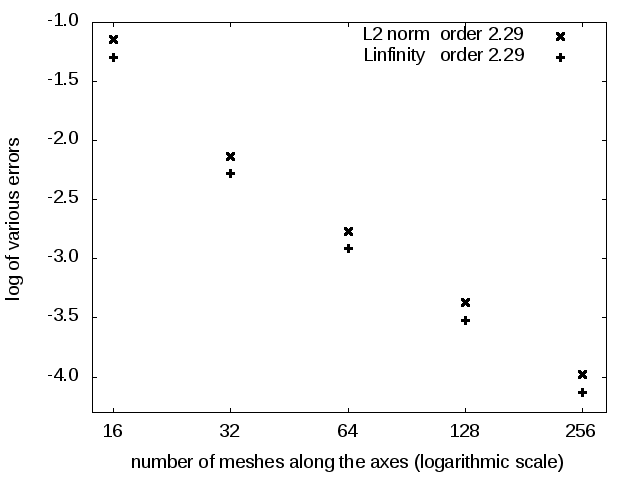}
\caption{  Two-dimensional wave with wave vector $ \, k = 2  \pi \, (1, \, 1 ) $, 
$ \, g \simeq  6  $, $ \,  2 \,  | k | \, c_0 \simeq 5.924 $. 
Convergence towards the damped acoustic model (\ref{damped-acoustic}).  }
  \label{manip02long-thermic} \end{figure}

\noindent
The initial condition is presented in Figure~\ref{d2q9periodic-densite-m021280}. 
The autocorrelation of density 
\moneqstar 
\displaystyle \Gamma (t) \equiv {{ \displaystyle  \int_\Omega \rho(x,\, t) \, \rho(x,\, 0) \, \dd x}\over
{ \displaystyle  \int_\Omega | \rho(x,\, 0) |^2 \, \dd x} }
\monendstar 
is typical  of a diffusion process as shown in 
Figure~\ref{manip02long-autocorrelation}. 
The convergence for simple dyadic meshes is presented in Figure~\ref{manip02long-thermic}.

   \vspace {-.5cm}   

   \vspace {-.5cm}   

   \vspace {-.5cm}   \begin{figure}    [H] \centering
\includegraphics [height=.55 \textwidth] {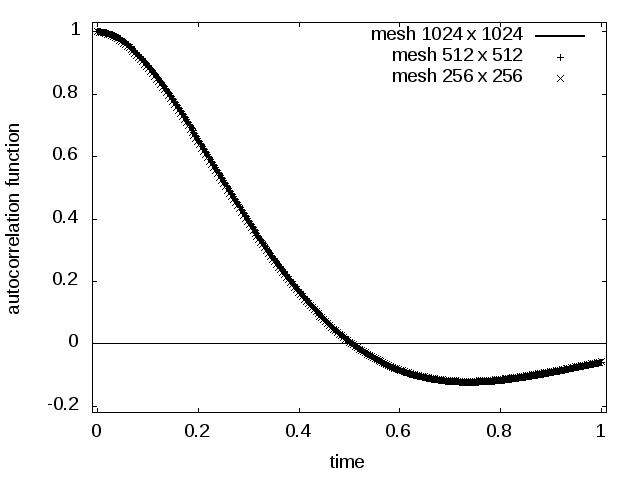}
\caption{  Two-dimensional wave with wave vector $ \, k = 2  \pi \, (1, \, 1 ) $, 
$ \, g \simeq  5.6470  $.  
Autocorrelation of density with  $ \,  2 \,  | k | \, c_0 \simeq 5.924 \, $
for various meshes.    }
  \label{manip0304-autocorrelation} \end{figure}

   \vspace {-.5cm}   \begin{figure}    [H] \centering
\includegraphics [height=.55 \textwidth] {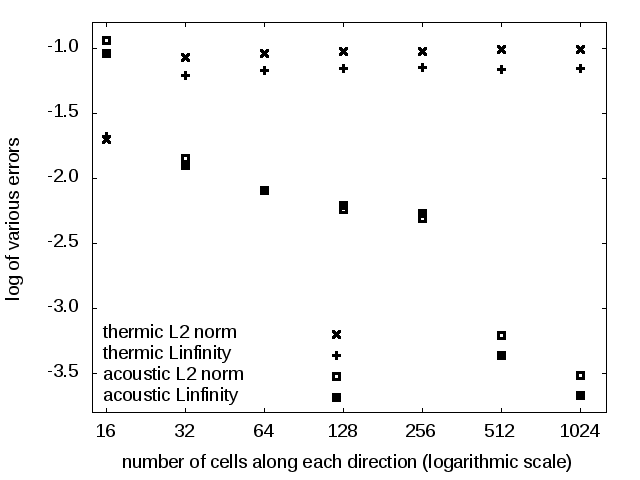}
\caption{  Two-dimensional wave with  vector $ \, k = 2  \pi \, (1, \, 1 ) $, 
$ \, g \simeq  5.6470  $.  
Convergence  towards the damped acoustic model (\ref{damped-acoustic}) 
 with  $ \,  2 \,  | k | \, c_0 \simeq 5.924 $. 
The order of convergence is 1.27 for the $ \, L^2 \, $  norm and 1.30 
in norm  $ \, L^\infty $.  }
  \label{manip0304-thermic-acoustic} \end{figure}

\newpage 
\noindent
A second numerical  experiment has been conducted. 
We keep the same domain  $ \, [ 0 ,\, 2 \, \pi ] ^2 \, $ with the same initial condition 
$ \, \rho = \cos \big( ( 2  \pi \, (x+y) \big) $. 
Then $ \, k = 2 \, \pi \, (1, \, 1 ) \,  $ and the right-hand  side of 
(\ref{propagative-acoustic})  
is equal to  $ \,  2 \,  | k | \, c_0 \simeq 5.924 $. 
We change the value 
of the diffusivity  $ \, \kappa \, $ introduced in (\ref{diffusivite})
to $ \, \kappa = {17\over288} \simeq 0.05903 $.   
We keep the traditional value $ \, c_0 = {{1}\over\sqrt{3}} $. 
Then the dissipation coefficient $ \, g \, $ (see (\ref{c0-gg})) is now    
 $ \, g \simeq 5.6470 $.     
Then the damped acoustic model  (\ref{damped-acoustic}) 
 exhibits a propagative mode in this case. 
The autocorrelation function is presented in Figure~\ref{manip0304-autocorrelation}. 
The convergence curve is depicted in Figure~\ref{manip0304-thermic-acoustic}.
We observe that this convergence is not regular. An extra-fine mesh
with dimensions  $1024\times 1024$ 
 has been  necessary in order to confirm the  order of accuracy.

%
%

\monitem 
 Complementary experiments for an initial Gaussian 

\noindent 
We have compared the scalar D2Q9 lattice Boltzmann scheme 
with  acoustic scaling with  numerical solutions of the heat equation 
(\ref{chaleur}) as presented  in Section~3 and with   HaWAY simulations
of the damped acoustic system (\ref{damped-acoustic})   in Section~4.   
We consider again the first geometry studied in this contribution, 
{\it id est} the  square $ \, \Omega = \, [ -1 , \, 1 ]^2 \, $ 
 with periodic boundary conditions. 
An initial Gaussian profile  (\ref{gaussienne}) is given at $ \, t = 0 $. 
Two numerical experiments have been considered: a quite viscous one 
with imposed diffusivity  $ \,  \kappa = 0.15 \, $ 
and another  one with $ \, \kappa =$~$ 0.015 $.  
The numerical parameters are displayed in Table~\ref{table-4}.

\begin{table}  [H]     \centering
 \centerline { \begin{tabular}{|c|c|c|c|c|c|c|c|c|}    \hline 
 number of cells & $ 13^2 $  & $ 27^2 $ & $ 55^2 $  &  $ 111^2 $ 
& $ 223^2 $  & $ 447^2 $  & $ 895^2 $ & $ 1791^2 $  \\   \hline 
$ s_J \, {\rm with} \, \kappa = 0.15 $  &   $  0.292  $ &  $  0.152 $ &  $ 0.0777 $  &  $ 0.0392 $ &
$ 0.0197 $  &  $ 0.00989 $ & &  \\   \hline 
$  s_J \, {\rm with} \, \kappa = 0.015 $  &   $  1.262 $ & $  0.903 $ &  0.575 &   0.333 &   0.181 
&  0.0947 &  0.0484 &  0.0245  \\   \hline  
\end{tabular} }  
\caption{Initial Gaussian.  D2Q9 numerical experiments with acoustic scaling ; values of $\, s_J \, $
for the viscosities  $ \, \kappa = 0.15 \, $ and  $ \, \kappa = 0.015 $. } \label{table-4} \end{table}


   \vspace {-.5cm}   \begin{figure}    [H] \centering
\includegraphics [ trim = {70mm 0mm 70mm 0mm}, clip, height=.305 \textwidth] {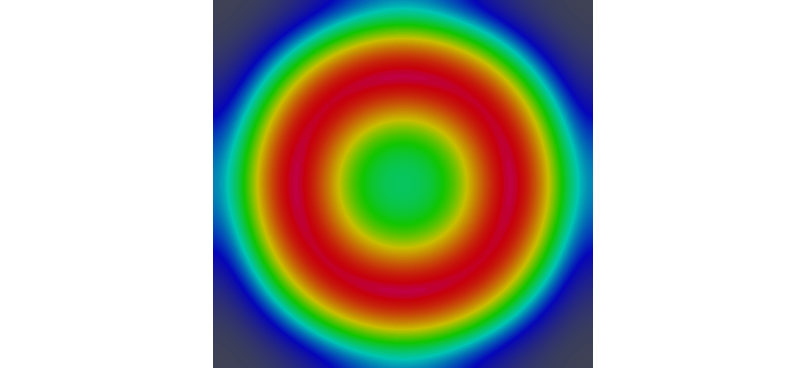}   
$ \!\!\!\!\!\!  $ 
\includegraphics [ trim = {70mm 0mm 70mm 0mm}, clip, height=.305 \textwidth] {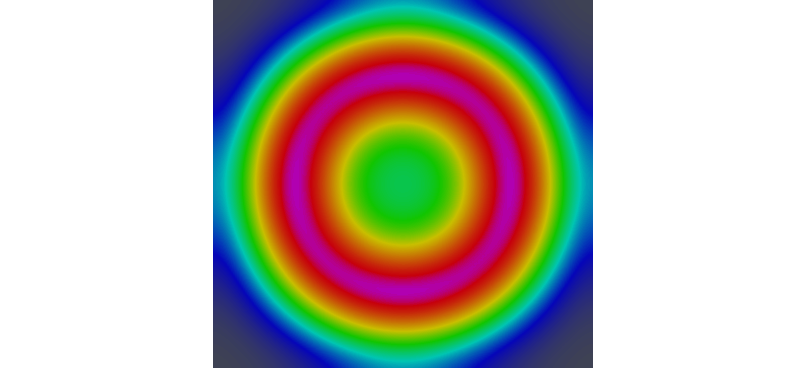}  
$ \!\!\!\!\!\!  $ 
\includegraphics [ trim = {70mm 0mm 70mm 0mm}, clip, height=.305 \textwidth] {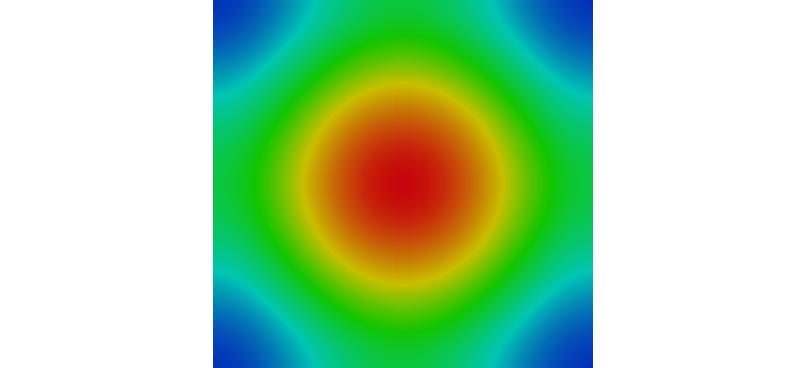} 
\caption{Initial Gaussian, $\, \kappa = 0.15 $. Three simulations done for the 
damped acoustic model with finite differences  HaWAY discretization (left), 
D2Q9 lattice Boltzmann scheme (middle),  and for the 
heat equation  with finite differences (right). Approximate solutions are presented for  time $ \, T = 1 \,$ 
with $\, 111 \times 111 \,$ meshes. }
  \label{111-111-gaussian-acoustic-t100} \end{figure}

   \vspace {-.5cm}   \begin{figure}    [H] \centering
\includegraphics [height=.55 \textwidth] {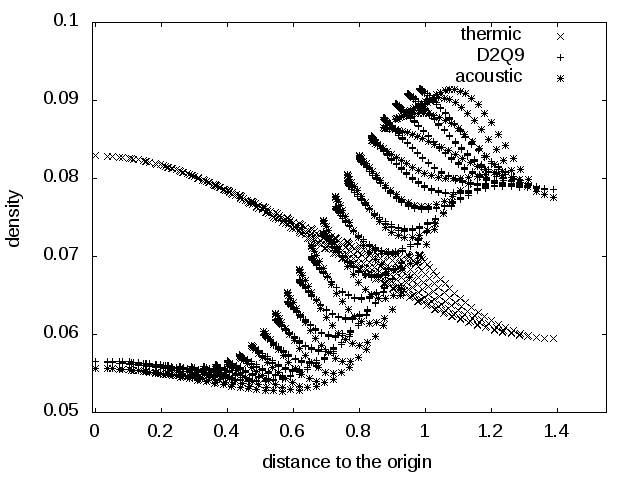}
\caption{Initial Gaussian, $\, \kappa = 0.15 $.  Three simulations done for the
damped acoustic model with finite differences  HaWAY discretization, 
D2Q9 lattice Boltzmann scheme  and 
heat equation  with finite differences. Density field at time = 2 of a 
 $\, 55 \times 55 \,$ mesh. }
  \label{manip5-time-2-rhor} \end{figure}

   \vspace {-.5cm}   \begin{figure}    [H] \centering
\includegraphics [height=.55 \textwidth] {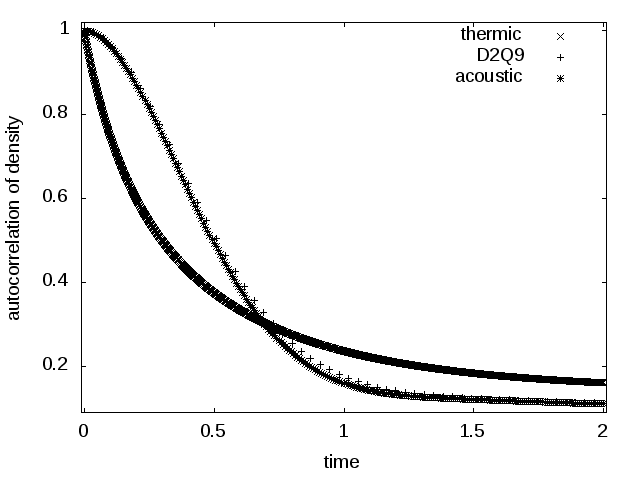}
\caption{Initial Gaussian, $\, \kappa = 0.15 $.  Three simulations done for the
damped acoustic model with finite differences  HaWAY discretization, 
D2Q9 lattice Boltzmann scheme  and 
heat equation  with finite differences.  Autocorrelation of density for a  
 $\, 55 \times 55 \,$ mesh. }
  \label{manip5-time-2-autocor} \end{figure}

   \vspace {-.5cm}   \begin{figure}    [H] \centering
\includegraphics [height=.55 \textwidth] {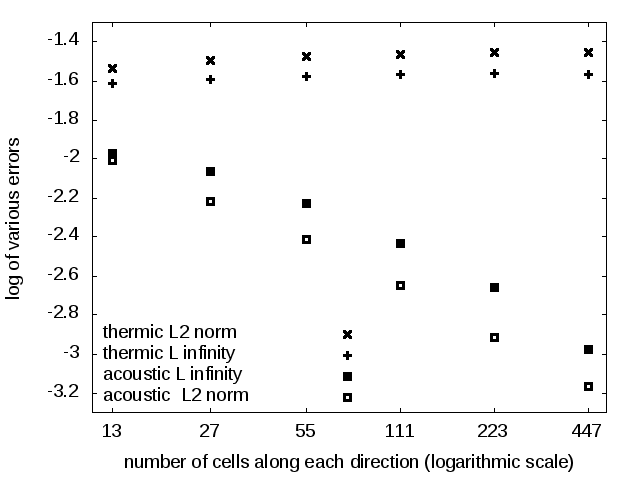}
\caption{Initial Gaussian, $\, \kappa = 0.15 $. Damped acoustic model with HaWAY finite differences, 
D2Q9 lattice Boltzmann scheme  and heat equation.
The order of convergence at time~$ = 2$ towards  damped acoustic is 0.756  for the $ \, L^2 \, $  norm and 0.653  
in norm  $ \, L^\infty $.  } 
  \label{manip5-time-2-convergence} \end{figure}

\noindent 
The results for the first test case with $ \, \kappa = 0.15 \, $ 
are presented in Figures \ref{111-111-gaussian-acoustic-t100}, 
\ref{manip5-time-2-rhor},   \ref{manip5-time-2-autocor} and~\ref{manip5-time-2-convergence}. 
In Fig.~\ref{111-111-gaussian-acoustic-t100},
a qualitative view of the numerical result on a given mesh shows that the scalar D2Q9 scheme   
and the HaWAY scheme for damped acoustic are closer to each other than they 
are to the solution of the heat equation. 
The three profiles of density are shown in  Fig.~\ref{manip5-time-2-rhor} and a comparison of autocorrelation 
functions in Fig.~\ref{manip5-time-2-autocor}.  Even on a relatively coarse mesh,   
the conclusion is the same and our new asymptotic analysis of  the acoustic system 
(\ref{damped-acoustic})  is consistent  with the numerical results. 
Last but not least, both  the error between D2Q9 and thermics on one hand, and that between 
D2Q9 and damped acoustics  on the other hand 
are  displayed in   Fig.~\ref{manip5-time-2-convergence}.   
The error between the lattice Boltzmann scheme and the damped acoustic results tends to zero 
whereas the error  between D2Q9 and the thermic model remains stationary. 

\begin{figure}    [H] \centering
\includegraphics [height=.55 \textwidth] {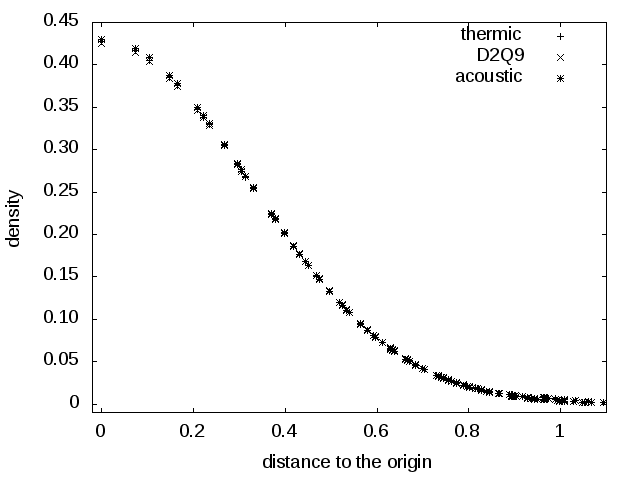}
\caption{Initial Gaussian, $\, \kappa = 0.015 $.  Three simulations done  for the
damped acoustic model with finite differences  HaWAY discretization, 
D2Q9 lattice Boltzmann scheme  and for the
heat equation  with finite differences.  Density field at time = 2 for a $\, 27 \times 27 \, $ mesh.}
  \label{manip7-time-2-rhor} \end{figure}

\begin{figure}    [H] \centering
\includegraphics [height=.55 \textwidth] {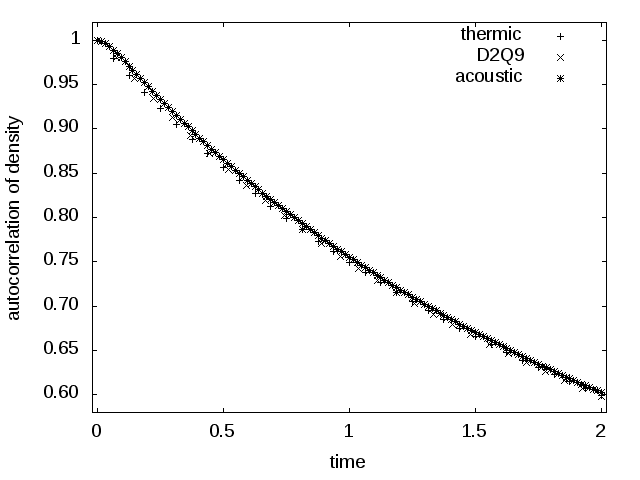}
\caption{Initial Gaussian, $\, \kappa = 0.015 $.  Three simulations done for the 
damped acoustic model with HaWAY finite differences, 
D2Q9 lattice Boltzmann scheme  and for the 
heat equation  with finite differences.  Autocorrelation of density at time = 2 for a  
 $\, 27 \times 27 \,$ mesh. }
  \label{manip7-time-2-autocor} \end{figure}

\begin{figure}    [H] \centering
\includegraphics [height=.55 \textwidth] {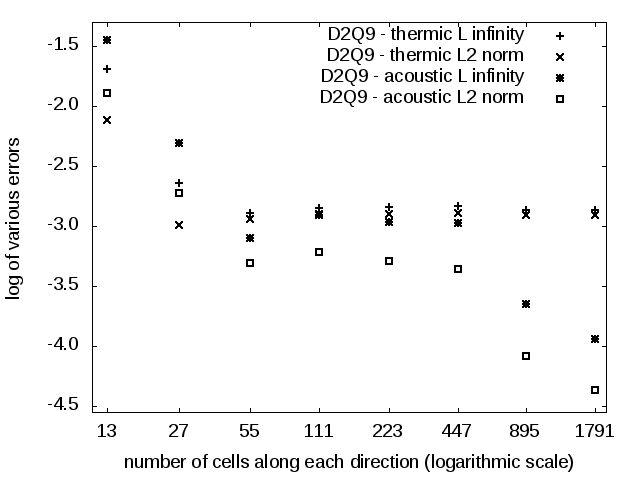}
\caption{Initial Gaussian, $\, \kappa = 0.015 $.  Three simulations done for the
damped acoustic model with finite differences  HaWAY discretization, 
D2Q9 lattice Boltzmann scheme  and for the
heat equation  with finite differences. 
The order of convergence at time~$ = 2$ towards  damped acoustic is 0.954 for the $ \, L^2 \, $  norm and 0.934  
in norm  $ \, L^\infty $.  } 
  \label{manip7-time-2-convergence} \end{figure}

\noindent 
The second numerical experiment with $ \, \kappa = 0.015 \, $ is presented 
in Figures~\ref{manip7-time-2-rhor}, 
\ref{manip7-time-2-autocor}, and  \ref{manip7-time-2-convergence}. At time $ \, =2 \,$ 
on a relatively  coarse mesh, 
the three numerical solutions can not be distinguished as shown in Fig.~\ref{manip7-time-2-rhor}. 
It is also the case 
for the autocorrelation function as presented in Fig.~\ref{manip7-time-2-autocor}. 
The numerical convergence is delicate for this test case. During one decade of mesh refinement, the three methods
present very close results as shown in Fig.~\ref{manip7-time-2-convergence}. 
Two additional computations on $\, 895 \times 895 \, $ and $ \, 1791 \times 1791 \, $ 
refined meshes have been necessary to demonstrate the 
convergence of the scalar D2Q9 scheme towards the damped acoustic system.  
Observe that the most refined mesh contains more than 3~millions cells!

   \vspace {-.5cm}   

   \vspace {-.5cm}   

   \vspace {-.7cm}   

\bigskip \bigskip   \noindent {\bf \large    6) \quad  Conclusion  }   

\noindent 
We have first considered the scalar D2Q9 lattice Boltzmann scheme with diffusive scaling.
Our experiments confirm numerical convergence to the solution of the heat equation. 
Of course, the mathematical proof of this numerical fact has now to be established. 

\noindent  
We have also studied convergence properties of the scalar D2Q9  scheme
with an acoustic scaling for the diffusion of a Gaussian profile, 
when it is supposed to approximate  diffusion problems.
Our numerical  experiments show  
consistent   
results  with the  diffusion equation solution    
 for  relaxation parameters that are not too small, $ s_J \geq 0.5 $ typically. 
When this relaxation coefficient is very small, however,  numerical convergence 
is defective  for the  diffusion of a Gaussian.

\noindent  
For very small values of the relaxation parameter, 
 the  asymptotic analysis has been revised  when the physical diffusion  is given.
We have developed a new analysis of the  lattice Boltzmann method
using  the dispersion equation and Gaussian elimination when  relaxation parameters 
can tend to zero. 
This asymptotic analysis shows that a  damped acoustic   model is
emergent  at first order. 
%
Complementary numerical experiments (see Figures~\ref{manip02long-thermic},
\ref{manip0304-thermic-acoustic}, \ref{manip5-time-2-convergence}
and~\ref{manip7-time-2-convergence})  show the numerical convergence 
of the  D2Q9 lattice Boltzmann scheme with acoustic scaling  
and a relaxation coefficient $\, s_J \,$ determined in such 
a way that the usual relation (\ref{mu}) is satisfied, 
towards the  damped acoustic system. 
Due to the mathematical convergence of  the lattice Boltzmann scheme
with diffusive scaling  \cite{JY15}, this result is unexpected, as pointed in the title. 

\noindent 
The results presented here can be interpreted physically in terms of frequency dependent transport coefficients 
that should be used when the time scale of the macroscopic phenomenon under study 
is not very large compared to microscopic time scales. 
%
Future study should focus on the extension of this analysis to second order. 
A natural extension of this question concerns  lattice Boltzmann models   
conserving {\it a priori} both mass and momentum. 
Our preliminary results show that a 
system of five partial differential equations is emergent 
in the case of  two space dimensions. 
This question will be studied in a forthcoming contribution.

%
\bigskip \bigskip   \noindent {\bf  \large  Acknowledgments }

\noindent 
The authors thank the Fondation Math\'ematique Jacques Hadamard for funding our collaboration.
This work is also partially supported  by the French ``Climb'' Oseo project.  
Last but not least, 
the authors thank the referees for their precise comments on the first draft of this contribution.

\bigskip \bigskip     \noindent {\bf \large    Appendix.  
HaWAY staggered finite differences } 

\noindent 
We consider the acoustic model proposed in Eq.~(\ref{damped-acoustic}): 
\moneq \label{damped-acoustic-2d} 
  \left \{ \begin {array}{l}
\displaystyle  {{\partial \rho}\over{\partial t}}  + {{\partial J^x}\over{\partial x}} 
+  {{\partial J^y}\over{\partial y}} = 0 
\\  \vspace{-4 mm} \\\displaystyle
 {{\partial J^x }\over{\partial t}}   +  c_0^2 \,\,  
{{\partial \rho}\over{\partial x}} + g  \, J^x = 0 
\\  \vspace{-4 mm} \\\displaystyle
 {{\partial J^y }\over{\partial t}}   +  c_0^2 \,\,  
{{\partial \rho}\over{\partial y}} + g  \, J^y = 0  \,. 
\end{array} \right. \monend
%

\vspace {-.5cm}  \begin{figure}   [H]     \centering
\includegraphics [width=.55 \textwidth, angle=0] {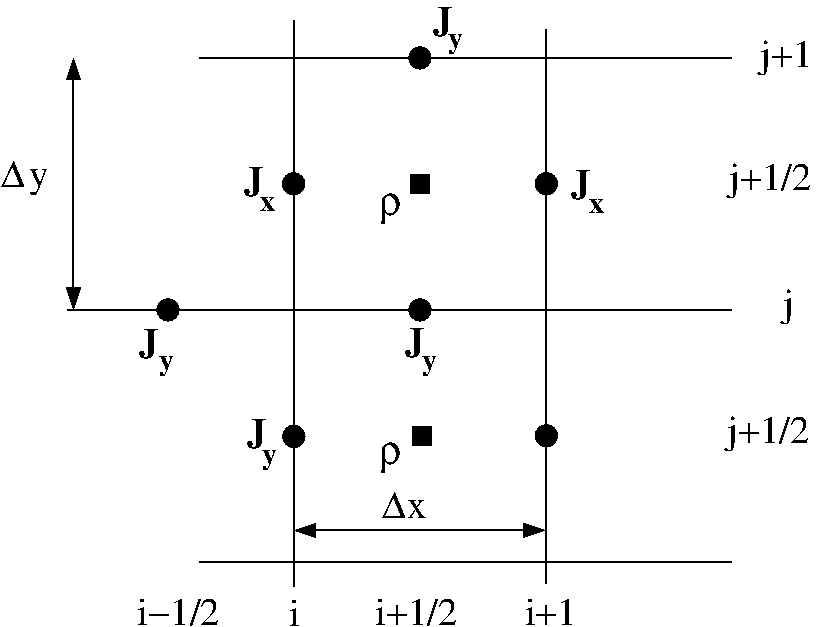}
\caption{HaWAY grid for staggered (density, momentum) finite differences   }
  \label{grille-HaWAY} \end{figure}

\noindent 
Given a spatial grid  $ \, \Delta x $, $ \, \Delta y \, $  and a time step $ \, \Delta t  $, 
we consider  integer multiples of these parameters for the discretization in  space and time. 
The density $ \,  \rho  \,$ is approximated at half-integer  vertices in space and integer points in time
whereas the momentum  $ \,  J^x  \,$  (respectively~$ \, J^y $)  
is approximated at integer nodes (respectively half-integer nodes) in the $x$-direction, 
 semi-integer nodes  (respectively integer nodes) in the $y$-direction, and  half-integer values in time:
\moneq \label{variables-discretes}  
 \rho \, \approx \, \rho^n_{i+1/2,\, j+1/2} \, , \quad J^x  \, \approx \, J^{x, \, n+1/2}_{i,\, j+1/2} 
\, , \quad J^y  \, \approx \, J^{y, \, n+1/2}_{i+1/2,\, j}  \, . 
\monend
The Figure~\ref{grille-HaWAY} gives an illustration of this classical choice 
\cite{Ar66, HaW65, Ye65}. 

\monitem 
We discretize the first equation of Eqs.~(\ref{damped-acoustic-2d}) with a four-point centered 
finite-difference schemes around the vertex 
$ \, \big( (i+{1\over2}) \, \Delta x \,,\, (j+{1\over2}) \, \Delta y 
 \,,\,  (n+{1\over2}) \, \Delta t \big) $: 
\moneq \label{masse-discrete}  
{{\rho^{n+1}_{i+1/2,\, j+1/2} - \rho^n_{i+1/2,\, j+1/2} }\over{\Delta t}}  + 
 {{ J^{x, \, n+1/2}_{i+1,\, j+1/2} - J^{x, \, n+1/2}_{i,\, j+1/2} }\over{\Delta x}}   + 
 {{ J^{y, \, n+1/2}_{i+1/2,\, j+1} - J^{y, \, n+1/2}_{i+1/2,\, j} }\over{\Delta y}}  \,=\, 0  \, . 
\monend
We use the same approach for the discretization of the second  equation of Eqs.~(\ref{damped-acoustic-2d}) 
 around the node  
$ \, \big( i \, \Delta x \,,\, (j+{1\over2}) \, \Delta y  \,,\,  n \, \Delta t \big) $: 
\moneq \label{impulsion-x-discrete}  \displaystyle 
{{  J^{x, \, n+1/2}_{i,\, j+1/2} -  J^{x, \, n-1/2}_{i,\, j+1/2} }\over{\Delta t}}  +  
{{c_0^2}\over{\Delta x}} \, \big( \rho^n_{i+1/2,\, j+1/2} - \rho^n_{i-1/2,\, j+1/2} \big)   +  
g \,  J^{x, \, n}_{i,\, j+1/2} \,=\, 0  
\monend
and the third  equation of Eqs.~(\ref{damped-acoustic-2d}) 
 around the node  
$ \, \big(  (i+{1\over2}) \, \Delta x \,,\, j \, \Delta y  \,,\,  n \, \Delta t \big) $: 
\moneq \label{impulsion-y-discrete}  \displaystyle 
{{  J^{y, \, n+1/2}_{i+1/2,\, j}  - J^{y, \, n-1/2}_{i+1/2,\, j} }\over{\Delta t}}  +  
{{c_0^2}\over{\Delta y}} \, \big( \rho^n_{i+1/2,\, j+1/2} - \rho^n_{i+1/2,\, j-1/2} \big)   +  
g \,   J^{y, \, n}_{i+1/2,\, j} \,=\, 0  \, . 
\monend
We interpolate the momentum at integer time vertices  with a simple average: 
\moneqstar 
 J^{x, \, n}_{i,\, j+1/2} =  {{1}\over{2}} \,  \big(    J^{x, \, n+1/2}_{i,\, j+1/2}
 +  J^{x, \, n-1/2}_{i,\, j+1/2}  \big)  \,, \,\, 
  J^{y, \, n}_{i+1/2,\, j}  =  {{1}\over{2}} \,  \big(   J^{y, \, n+1/2}_{i+1/2,\, j}  + J^{y, \, n-1/2}_{i+1/2,\, j} \big)  
\, .  
\monendstar 
We incorporate these expressions into  the relations Eq.~(\ref{impulsion-x-discrete}) 
and Eq.~(\ref{impulsion-y-discrete}). We obtain 
\moneq \label{impulsion-x-discrete-bis}  
\Big(  {{1}\over{\Delta t}} + {{g}\over{2}} \, \Big) \,   J^{x, \, n+1/2}_{i,\, j+1/2} 
+ {{c_0^2}\over{\Delta x}} \, \big( \rho^n_{i+1/2,\, j+1/2} - \rho^n_{i-1/2,\, j+1/2} \big) = 
\Big(  {{1}\over{\Delta t}} - {{g}\over{2}} \, \Big) \,  J^{x, \, n-1/2}_{i,\, j+1/2} \, 
\monend
and 
\moneq \label{impulsion-y-discrete-bis}  
\Big(  {{1}\over{\Delta t}} + {{g}\over{2}} \, \Big) \,   J^{y, \, n+1/2}_{i+1/2,\, j}  + 
{{c_0^2}\over{\Delta y}} \, \big( \rho^n_{i+1/2,\, j+1/2} - \rho^n_{i+1/2,\, j-1/2} \big) = 
\Big(  {{1}\over{\Delta t}} - {{g}\over{2}} \, \Big) \, J^{y, \, n-1/2}_{i+1/2,\, j} \, . 
\monend

The numerical scheme is now entirely defined for internal nodes. 
In this study we  have used periodic boundary conditions. 


\bigskip \bigskip      \noindent {\bf  \large  References }   


\end{document}